\documentclass[11pt,twoside]{article}
\usepackage{latexsym, amssymb, amsmath, theorem}
\usepackage{epsfig}
\numberwithin{equation}{section}

\theoremstyle{change}
\theorembodyfont{\itshape}
\newtheorem{theorem}{Theorem}[section]
\newtheorem{proposition}[theorem]{Proposition}
\newtheorem{lemma}[theorem]{Lemma}
\newtheorem{corollary}[theorem]{Corollary}

\theorembodyfont{\rmfamily}
\newtheorem{definition}[theorem]{Definition}

\newtheorem{remark}[theorem]{Remark}
\newenvironment{proof}{{\noindent \textbf{Proof}\,\,}}{\hspace*{\fill}$\Box$\medskip}

\title{Restricted version of the infinitesimal
Hilbert 16th problem\protect\footnote{AMS Classification 
58F21 (14K20 34C05)}}

\author{A. Glutsyuk\protect\footnote{Laboratoire
J.-V.Poncelet (UMI 2615 du CNRS et
l'Universit\'e Ind\'ependante de Moscou). Permanent address:
CNRS, Unit\'e de Math\'ematiques Pures et Appliqu\'ees,
M.R., \'Ecole Normale Sup\'erieure de Lyon, 46 all\'ee d'Italie,
69364 Lyon 07, France.},
Yu. Ilyashenko\protect\footnote{Moscow State and Independent Universities, Steklov Math.
Institute, Moscow; Cornell University, US.}
\protect\footnote{Both authors were supported by part by RFBR grant
02-02-00482.
The second author was supported by the grant NSF 0400495.}
}

\begin{document}
\maketitle

\def\Cal{\mathcal}
\def\Bbb{\mathbb}
\def\e {\varepsilon }
\def\nbd {neighborhood }
\def\nbds {neighborhoods }
\def\ph {\varphi }
\def\a {\alpha }
\def\de {\delta }
\def\g {\gamma }
\newcommand\diffeo {diffeomorphism }
\newcommand\ihp {Infinitesimal Hilbert 16th Problem }
\newcommand\un {uniformization }
\renewcommand\inf {infinitesimal }
\newcommand\uni {uniformization }
\newcommand\hn {$\mathcal H_{n}$}
\newcommand\bui{B_{K,U} (I)}

\tableofcontents

\section{Introduction}

\subsection{Restricted \ihp}

The original \ihp is stated as follows. Consider a real polynomial
$H$ in two variables of degree $n+1.$ The space of all such
polynomials is denoted by \hn.

Connected components of closed level curves of $H$ are called {\it
ovals} of $H.$ Ovals form continuous families, see Fig. 1. Fix one
family of ovals, say $\Gamma,$ and denote by $\gamma_t$ an oval of
this family
that belongs to the level curve $\{H=t\}.$

\begin{figure}[ht]
  \begin{center}
   \epsfig{file=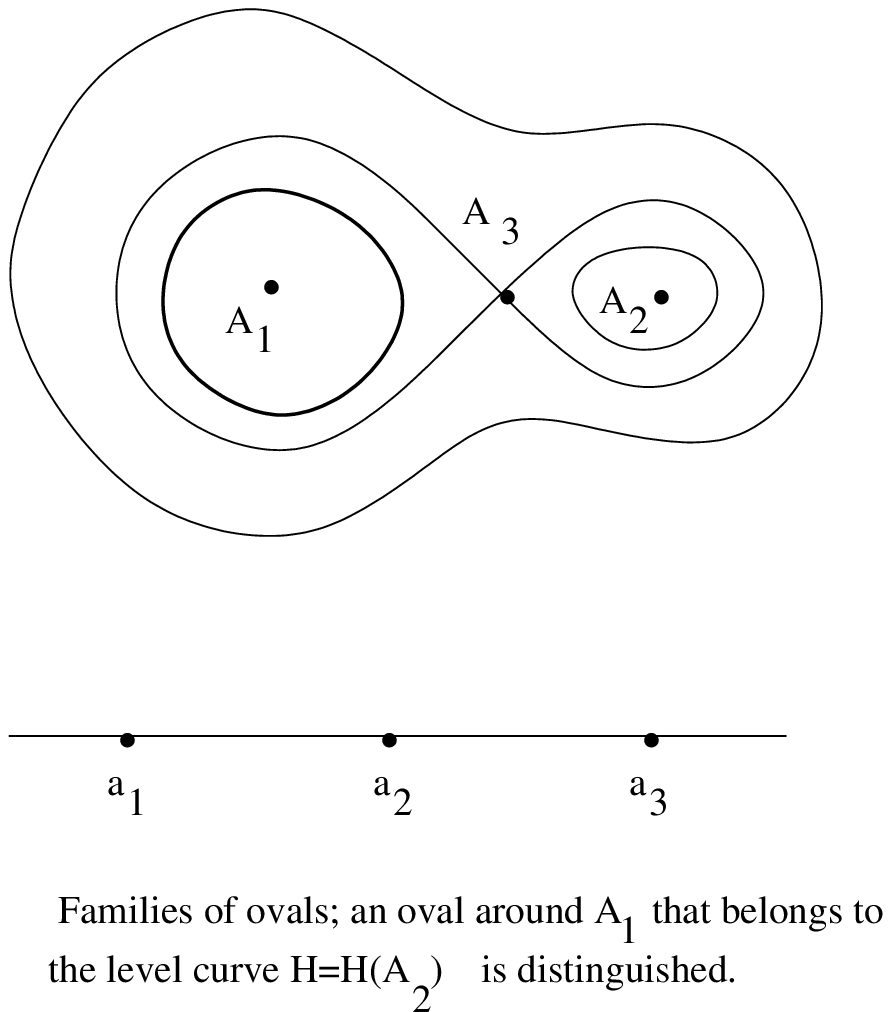}
    \caption{}
    \label{fig:1}
  \end{center}
\end{figure}

Consider a polynomial one-form
$$ \omega = A dx + B dy $$
with polynomial coefficients of degree at most $n.$ The set of all
such forms is denoted by $\Omega_n^.$ The main object to study
below is the integral
\begin{equation}
I(t) = \int_{\gamma_t} \omega.        \label{1.1}
\end{equation}

\medskip
{\bf \ihp.} {\it Let $H$ and $\omega $ be as above. Find an upper
bound of the number of isolated real zeros of integral (\ref{1.1}) for a
polynomial $H \in \mathcal H_{n}$ and any family $\Gamma $ of real
ovals of $H.$ The estimate should be uniform in $\omega $ and $H,$
thus depending on $n$ only.}

\medskip

This problem stated more than 30 years ago is not yet solved.
The existence of such a bound was proved by A.N.Varchenko \cite{33}
and A.G.Khovanskii \cite{15}. A
weaker version of the problem is called {\it restricted.} In order
to formulate it we need the following

\begin{definition} \label{d1.1} A polynomial $H \in \mathcal H_{n}$ is {\it
ultra-Morse} provided that it has $n^2$ complex Morse critical
points with pairwise distinct critical values, and the sum $h$ of
its higher order terms has no multiple linear factors.
\end{definition}


Denote by $\mathcal U_n$ the set of all ultra-Morse polynomials in
\hn. The complement to this set is denoted by $\Sigma_n $ and called
{\it the discriminant set}. The integral (\ref{1.1}) may be identically
zero. The following theorem shows that for ultra-Morse polynomials
this may happen by a trivial reason only.

\begin{theorem} ({\bf Exactness theorem \cite {11, 12, 30a}}).

Let $H$ be a real ultra-Morse polynomial of degree higher than 2.
Let the integral (\ref{1.1}) be identically zero for some family
of real ovals of the polynomial $H.$ Then the form $\omega$ is
exact: $ \omega = df.$
\end{theorem}
Denote by $\Omega_n^*$ the set of all non-exact polynomial
one-forms from $\Omega_n.$

\medskip
{\bf Restricted version of the \ihp.} {\it For any compact set
$\mathcal K \subset U_{n}$ find an upper bound of the number of all
real zeros of the integral (\ref{1.1}) over the ovals of the polynomial
$H \in \mathcal K.$ The bound should be uniform with respect to $H \in
\mathcal K$ and $\omega \in \Omega_n^*.$ It may depend on $n$ and
$\mathcal K$ only.}
\medskip

This problem is solved in the present paper, and the explicit
upper bound is given in the next subsection.

\subsection{Main results}

To measure a gap between a compact set $\mathcal K \subset U_{n}$ and
the discriminant set $\Sigma_n$, let us first normalize
ultra-Morse polynomials by an affine transformation in the target
space. This transformation does not change the ovals of $H$, thus
the number of zeros of the integral (\ref{1.1}) remains unchanged.

Say that two polynomials $G$ and $H$ are equivalent iff
$$
G = a H + b, \ \ a > 0, \ b \in \Bbb C.
$$

\begin{definition} \label{d1.2} A polynomial is {\it balanced} if all its
complex critical values belong to a disk of radius $2$ centered at
zero, and there is no smaller disk that contains all the critical
values.
\end{definition}

\begin{remark} Any polynomial with at least two distinct
critical values is equivalent to one and unique balanced
polynomial. If the initial polynomial has real coefficients, then
so does the corresponding balanced polynomial.
\end{remark}

Define two positive functions on $\mathcal U_n$ such that at least one
of them tends to zero as $H$ tends to $\Sigma_n.$ For any compact
set $\mathcal K \subset \mathcal U_n$ the minimal values of these
functions on $\mathcal K$ form a vector in $\Bbb R^+ \times \Bbb R^+$
that is taken as a size of the gap between $\mathcal K$ and
$\Sigma_n.$

\begin{definition} \label{d1.3} For any $H \in \mathcal U_n$ let $c_1(H)$ be
$n$ multiplied by the smallest distance between two lines in the
locus of $h,$ the higher order form of $H.$ The distance between
two lines is taken in sense of Fubini-Study metric on the
projective line $\Bbb CP^1.$ Let $c'(H) = \min (c_1(H),1).$
\end{definition}
 Denote by $\mathcal V_n$ the set of all polynomials with more
than one critical value and more than one line in the locus of the
higher order homogeneous form. By Definition 1.1, $\mathcal U_n
\subset \mathcal V_n.$

\begin{definition} \label{d1.4} For any $H \in \mathcal V_n,$ let $G$ be the
balanced polynomial equivalent to $H.$ Let $c_2(H)$ be the minimal
distance between two critical values of $G$ multiplied by $n^2.$
Let $c''(H) = \min (c_2 (H),1).$
\end{definition}

Note that inequality $c'(H)c''(H) > 0$ is equivalent to the
statement that $H$ is ultra-Morse.

In what follows, we deal with balanced ultra-Morse polynomials
only. This may be done without loss of generality: any ultra-Morse
polynomial is equivalent to a balanced one; equivalent polynomials
have the same number of zeros of the integral (\ref{1.1}) over the same
family of ovals.

\medskip

{\bf Theorem A.} {\it Let $H$ be a real ultra-Morse polynomial of
degree $n + 1.$ Let $\Gamma = \{ \gamma_t \} $ be an arbitrary
continuous family of real ovals of $H.$ There exists a universal
positive $c$ such that the integral (\ref{1.1}) has at most $(1 - \log
c'(H))e^{\frac {c}{c''(H)}n^4}$ isolated zeros.}
\medskip

{\bf Appendix.} The statement of Theorem A holds with $c =
5.000.$
\medskip

An approach to the \ihp itself presented below motivates the
following complex counterpart of Theorem A, namely, Theorem B that
gives an estimate of the number of zeros of the integral (\ref{1.1}) in
the complex domain. Consider an ultra-Morse polynomial $H$ and let
\begin{equation}
\nu = \nu (H): = \frac {c''(H)}{4n^2}   \label{1.2}
\end{equation}
Fix any real noncritical value $t_0$ of $H,$
$$|t_0|<3,$$
whose distance to the complex critical values of $H$ is no less
than $\nu .$ Consider a real oval $\gamma_0 \subset \{ H= t_0\} .$
We suppose that such an oval exists. Let $a = a(t_0) < t_0 <
b(t_0) = b$ (or $a(H,t_0), \ b(H,t_0)$ for variable $H$) be the
nearest real critical values of $H$ to the left and to the right
from $t_0$ respectively; or $-\infty , +\infty $ if there are
none. Denote by $\sigma (t_0)$ the interval $(a(t_0), b(t_0))$ and
let $\Gamma (\gamma_0)$ be the continuous family of ovals that
contains $\gamma_0:$
\begin{equation}
\Gamma (\gamma_0) = \{ \gamma (t)\ | t \in \sigma (t_0), \ \gamma
(t_0) = \gamma_0\} .   \label{1.2a}
\end{equation}
The following cases for $(a,b) = \sigma (t_0)$ are possible:
$$
(a,b), \ b > a; \ (a, +\infty ); \ (-\infty , b).
$$
If $\lim {\text { top }}_{t\to a}\gamma (t)$ contains a critical
point of $H,$ then $a$ is a logarithmic branch point of $I.$ If
not, $a$ is called an apparent singularity. The same for $b.$ Let
$$
W =W(t_0,H)
$$
be the universal cover over the set of noncritical values of $H$
with the base point $t_0$ and the projection $\pi : W \to
\Bbb C.$ Let $D(t,r)$ be the disk centered at $t$ of radius $r.$
Denote by $a + re^{i\ph } \in W$ a point represented by a curve
$\Gamma_1 \Gamma_2 \subset W,$ where $\Gamma_1$ is an oriented
segment from $t_0$ to $t_1 = a + r \in \sigma (t_0), \ \Gamma_2 =
\{ a + re^{i\theta}\mid \theta \in [0,\ph ]\} ;$ $\Gamma_2$ is
oriented from $t_1$ to $t.$ In the same way $b - re^{i\ph } \in W$
is defined. Let
\begin{equation}
\Pi(a) = \{ a + re^{i\ph } \in W \mid 0 < r \le \nu , |\ph | \le
2\pi \} , \text { for } a \ne -\infty    \label{1.2b}
\end{equation}
$$
\Pi(b) = \{ b - re^{i\ph } \in W \mid 0 < r \le \nu , |\ph | \le
2\pi \} , \text { for } b \ne +\infty
$$
Let
$$
D(l,a) = \{ a + re^{i\ph } \in W \mid a + re^{\frac {i\ph }{l}}
\in \Pi (a)\}
$$
$$
D(l,b) = \{ b - re^{i\ph } \in W \mid b - re^{\frac {i\ph }{l}}
\in \Pi (b)\}
$$
Let $DP_R = DP_R(H,t_0)$ be the disk of
radius $R$ in the Poincar\' e metric of $W$ centered at $t_0.$
Denote
$$S_t=\{ H=t\}\subset\mathbb C^2.$$

For any real polynomial $H,$ the choice of a cycle $\gamma_0$
determines a family of ovals (\ref{1.2a}) over which the
integral (\ref{1.1})
is taken. When we want to specify this choice we write
$I_{H,\gamma_0}$ or $I_H$ instead of $I.$ The integral
$I_{H,\gamma_0}$ may be analytically extended not only as a
function of $t,$ but also as a function of $H.$ We stress that the
integral $I_{H,\gamma_0}(t)$ is taken over the oval
$\gamma_{H,t}(t) \subset S_t;$ the family of ovals
depends continuously on $H$ and $t$ and may be expended to complex
values of $t$ as elements of the homology group $H_1(S_t,
\mathbb Z)$ depending continuously on $H$ and $t.$

An analytic extension of the integral $I$ to $W$ is denoted by the
same symbol $I.$ For any positive $R$ and natural $l$ denote by $G
= G(l,R,H,t_0)$ the domain
$$
G = DP_R(H,t_0)\ \cup\ D(l,a(H,t_0))\ \cup\ D(l,b(H,t_0)).
$$

\medskip
{\bf Theorem B.} {\it For any real ultra-Morse polynomial $H,$ any
 real oval $\gamma_0$ of $H,$ any natural $l$ and any positive $R
> \frac {288n^4}{c''(H)},$ the number of zeros of the integral $I_{H,\gamma_0}$ in
$G = G(l,R,H,t_0),$  where $ t_0 = H\mid \gamma_0,$ is estimated
as follows:}
\begin{equation}
\# \{ t \in G(l,R,H,t_0)|I_{H,\gamma_0}(t) = 0\} \le (1 - \log
c'(H)) \cdot \big( e^{7R} + A ^{4800}e^{\frac {481l}{c''(H)}}\big)
, \ A = e^{\frac {n^4}{c''(H)}}. \label{1.3a}
\end{equation}
\medskip

The lower bound on $R$ in the statement of the theorem is
motivated by the remark in Subsection 1.7 below.

\subsection{An approach to a solution of the \ihp}

{\bf Conjecture.} {\it For any $n$ there exist $\delta (n), l(n),
R(n)$ with the following property. Let $H_0$ be an arbitrary real
polynomial from $\mathcal H_n, \ t_0$ be its real noncritical
value and $\gamma_0$ be a real oval of $H_0$ that belongs to $\{
H_0 = t_0\} $ (we suppose that such an oval exists). Let $I_H$ be
the integral (\ref{1.1}). The integral $I_H$ depends on $H$ as a
parameter. Let $t_1 \in \sigma (t_0), I_{H_0}(t_1) = 0$ and $t(H)$
be a germ of an analytic function defined by the equation
$I_H(t(H)) \equiv 0, \ t(H_0) = t_1.$ The required property is the
following. There exists a path $\lambda \subset \mathcal H_n$
depending on $H_0$ only starting at $H_0$ and ending at some $H_1
\in \mathcal H_n$ such that:
$$
c'(H_1) \ge \delta (n), \ c''(H_1) \ge \delta (n);
$$
$t_0$ is a noncritical value of all the polynomials along the path
$\lambda ;$

the analytic extension $t(H_1)$ of the function $t(H)$ along
$\lambda $ starting at the value $t_1$ belongs to the domain
$G(l(n), R(n), H_1, t_0).$}

\medskip

The conjecture above implies the solution of the Infinitesimal
16th Problem. Indeed, suppose that the conjecture is true. Let
$$
N(n) = (1 - \log \delta (n))e^{7R(n)+\frac {c_1n^4}{\delta
(n)}+\frac {c_2l(n)}{\delta (n)}}, c_1 = 4800, c_2 = 481
$$
Then the number of real zeros of integral $I_{H_0}$ can not exceed
$N(n).$ If not, any of real zeros of $I_{H_0}$ would be extended
along $\lambda $ up to a zero of a polynomial $H_1$ located in $G
= G(l(n), R(n), H_1,t_0).$ Thus the  number of zeros of the
integral $I_{H_1}$ in $G$ will exceed $N(n).$ But Theorem B
implies:
$$
\# \{ t \in G(l(n), R(n), H_1, t_0)\mid I_{H_1,t_0}(t) = 0\} \le
N(n),
$$
a contradiction.

\subsection{Historical remarks}

A survey of the history of \ihp may be found in \cite { cent},
and we will not repeat it here. In particular, a much weaker
version of Theorem A is claimed there as Theorem 7.7. The first
solution to restricted Hilbert problem was suggested in \cite
{22}. An explicit upper bound for the same numbers of zeros as in
Theorem A was suggested there as a tower of four exponents with
coefficients ``that may be explicitly written following the
proposed constructive solution.''  It is unclear how much efforts
is needed to write these constants down. Moreover, exponential of
a polynomial presented in Theorem A is much simpler (though still
very excessive) than the tower of four exponentials.

The result of \cite {22} is a crown of a series of papers
\cite{19} - \cite{21}. Solution to the restricted version of the \ihp
presented there is only one application of a vast theory.

This theory presents an upper bound of the number of zeros of
solutions to linear systems of differential equations. Similar
results for components of vector solutions to linear systems are
obtained. Abelian integrals are considered as solutions to
Picard-Fuchs equations.

On the contrary, our presentation is focused on the study of
Abelian integrals given by formula (\ref{1.1}) ``as they are'' and not
as solutions of differential equations.

\subsection{Quantitative algebraic geometry}

Our main tool is Growth-and-Zeros theorem for holomorphic
functions stated in the next subsection. It requires, in
particular, an upper bound of the integral under consideration. We
fix an integrand, say $w = x^ky^{n-k}dx.$ Depending on a scale in
$\Bbb C^2,$ a cycle $\gamma $ in the integral $\int_\gamma \omega
$ may be located in a small or in a large ball. According to this,
the integrand will be small or large. We want to estimate the
integral at a certain point of the universal cover $W$ represented
by an arc  that connects a base point $t_0$ with some point, say
$t,$ with $|t| \le 3.$ To make this restriction meaningful, the
scale in the range of the polynomial should be chosen; in other
words,  the polynomial should be balanced. The argument above
shows that it should be also {\it rescaled} in sense of the
following definitions.
\begin{definition} \label{d1.6} The
{\it norm} of a homogeneous polynomial is
the maximal value of its module on the unit sphere; this norm
is denoted by $\| h \|_{\max}.$
\end{definition}
\begin{definition} \label{d1.7} A balanced polynomial $H \in \Bbb
C[x,y]$ is {\it rescaled} provided that the norm of its higher order
form $h$ equals one: $||h||_{\max} = 1,$ and the origin is a
critical point for $H.$ Briefly, a balanced
rescaled polynomial will be called {\it normalized}.
\end{definition}
\begin{remark} Any ultra-Morse polynomial may be transformed to
a normalized one by affine transformations in the source
and target spaces (not in the unique way). The functions $c'$ and
$c''$ remain unchanged under such transformations.
\end{remark}
\begin{definition} \label{d1.8} We say that
the topology of a level curve $S_t =
H^{-1}(t)$ of a polynomial $H \in \mathcal H_n$ is {\it located in a
bidisk}
$$
D_{X,Y} = \{ (x,y) \in \Bbb C^2| |x| \le X, \ |y| \le Y\}
$$
provided that the difference $S_t \setminus D_{X,Y}$ consists of
$n +1= \deg H$ punctured topological disks, and the restriction of
the projection $(x,y) \mapsto x$ to any of these disks is a
biholomorphic map onto $\{ x \in \Bbb C| X < |x| < \infty \} .$
\end{definition}
\medskip

\def\baru{\overline U}

{\bf Theorem C \cite{G*}.}
{\it For a normalized polynomial, the
Hermitian basis in $\Bbb C^2$ may be so chosen that the  topology
of all level curves $S_t$ for $|t| \le 5$ will be located in a
bidisk $D_{X,Y}$ with}
$$
X\leq Y \le  {(c'(H))}^{-14n^3}n^{65n^3}=R_0.
$$
\medskip

This theorem is of independent interest, providing one of the
first results in {\it quantitative algebraic geometry.} On the
other hand, it implies upper estimates of Abelian integrals used
in the proof of Theorem A and required by the Growth-and-Zeros
theorem below.

In what follows, we describe the main ideas of the proof of a
simplified version of Theorem A, namely Theorem A1 stated below.
It provides an upper bound for the number of zeros of the integral
(\ref{1.1}) on a real segment that is $\nu $-distant from critical
values of $H$ and belongs to the disk
$\overline D_3 = \{ t\mid |t| \le 3\} ,$
thus being distant from infinity; recall that $\nu = \nu (H)$ is
given by (\ref{1.2}).

By the use of Theorem A1, we get in Section 4 an estimate of the
number of zeros of the integral $I_{H, \gamma_0}$ near the endpoints
of $\sigma (t_0),$ as well as near infinity (Theorem A2 stated
in 1.8). Together
with Theorem A1, this completes the proof of Theorem A. Theorem B
is split into two parts. The first one (Theorem B1 stated in 3.1)
is proved by
extending an upper bound given with the help of Theorem C from a
disk $|t| \le 5$ into a larger domain. The second one, Theorem B2
stated in 4.6,
is proved (in the same place) by the same tools as Theorem A2 that
include Petrov method and a so called KRY theorem. The latter one
is a recent result in one-dimensional complex analysis \cite{16, 31}.
Its improved version is proved by the
second author (Yu.S.Ilyashenko) in a separate paper \cite{KRY}
and stated in Section 4. In this form it provides a mighty
tool to estimate the number of zeros of analytic functions near
logarithmic singularities.

\subsection{Growth-and-Zeros Theorem for Riemann surfaces}

The idea of the proof of Theorem $A1$ is to consider an analytic
extension of the integral (\ref{1.1}) to the complex domain and to make
use of the following Growth-and-Zeros theorem.  The definition of
the intrinsic diameter used in the statement of the theorem is
recalled below. We need the following

\begin{definition} \label{d1.9}
Let $W$ be a Riemann surface, $\pi : W \to \Bbb
C$ be a holomorphic function (called projection) with non-zero
derivative. Let $\rho$ be the metric on $W$ lifted from $\Bbb C$
by projection $\pi .$ Let $U \subset W$ be a connected domain, and
$K \subset U $ be a compact set. For any $p \in U$ let $\e (p,
\partial U)$ be the supremum of radii of disks centered at $p,$
located in $U$ and such that $\pi $ is bijective on these disks.
The $\pi $-gap between $K$ and $\partial U,$ is defined as
$$
\pi \text {-gap }(K, \partial U) = \min_{p\in K} \e (p,\partial
U).
$$
\end{definition}

{\bf Growth-and-zeros theorem.} {\it Let $W, \pi , \rho$ be the
same as in Definition \ref{d1.9}. Let $U \subset W$ be a domain
conformally equivalent to a disk. Let $K \subset U$ be a path
connected compact subset of $U$ (different from a single point).
Suppose that the following two assumptions hold:

Diameter condition:
$$
{\text { diam }}_{int} K \le D;
$$

Gap condition:
$$
\pi \text {-gap} (K, \partial U) \le \e .
$$

\def\baru{\overline U}

Let $I$ be a bounded holomorphic function on $\bar U.$ Then}
\begin{equation}
\# \{ z \in K|I(z) = 0\} \le e^{\frac {2D}{\e }}\log\frac
{\max_{\baru}|I|}{\max_K|I|}   \label{1.3b}
\end{equation}
\medskip

The definition  of the intrinsic diameter is well known; yet we
recall it for the sake of completeness.
\begin{definition} \label{d1.10} The {\it intrinsic distance} between two
points of a path connected set in a metric space is the infinum of
the length of paths in $K$ that connect these points (if exists).
The {\it intrinsic diameter} of $K$ is the supremum of intrinsic
distances between two points taken over all the pairs of points in
$K.$
\end{definition}
\begin{definition} \label{d1.11} The second factor in
the right hand side
of (\ref{1.3b}) is called the Bernstein index of $I$ with respect to $U$
and $K$ and denoted $B_{K,U}(I):$
\begin{equation}
B_{K,U}(I) = \log \frac {M}{m}, \ M = \sup_U|I|, \ m = \max_K|I|.
\label{1.4}
\end{equation}
\end{definition}

\begin{proof} {\bf of the Growth-and-Zeros theorem.} The above theorem is
proved in \cite {13} for the case when $W = \Bbb C, \pi = Id.$ In
fact, in \cite {13} another version of (\ref{1.3b}) is proved with
(\ref{1.3b}) replaced by
\begin{equation}
\# \{ z \in K|I(z) = 0\} \le B_{K,U}(I)e^\rho ,  \label{1.5n}
\end{equation}
where $\rho $ is the diameter of $K$ in the Poincar\' e metric of
$U.$ In this case it does not matter whether $U$ belongs to $\Bbb
C$ or to a Riemann surface.

\begin{proposition} \label{p1.1}
Let $K, U$ be two sets in the Riemann
surface $W$ from Definition \ref{d1.9}, and let the Diameter and Gap
conditions from the Growth-and-Zeros theorem hold. Then the
diameter of $K$ in the Poincar\' e metric of $U$ admits the
following upper estimate:

\begin{equation}
\rho \le 2D/\e .   \label{1.5a}
\end{equation}
\end{proposition}
\begin{proof}  By the monotonicity property of the Poincar\' e
metric, the length of any vector $v$ attached at any point $p \in
K$ is no greater than two times the Euclidean length of $v$
divided by the $\pi $-gap between $K$ and $\partial U.$ This
implies (\ref{1.5a})
\end{proof}

Together with (\ref{1.5n}), this proves (\ref{1.3b}).
\end{proof}

\subsection{Theorem A1 and Main lemma}

In what follows, $H$ will be an ultra-Morse polynomial unless the
converse stated. Consider a normalized
polynomial $H$. Let $a_j$ be its complex critical values,
$j=1,\dots,n^2; \ \nu, \ t_0, \ W$ and $\pi$ be the same as in
1.2. Let $I$ be the  integral (\ref{1.1}) as in Theorem A (well defined
for $t=t_0$). It admits an analytic extension to $W$, which will
be denoted by the same symbol $I.$

Let $a = a (t_0), b = b (t_0)$ be the same as in 1.2, and $\nu $
be from (\ref{1.2}). Let
$$
l(t_0) = \begin{cases} a+\nu \text { for } a \not = -\infty \\
         -3 \text { for } a = -\infty, \end{cases}
$$
$$
r(t_0) = \begin{cases} b-\nu \text { for } b \not = +\infty \\
         3 \text { for } b = +\infty  . \end{cases}
$$
Let
$$
\sigma (t_0, \nu ) = [l(t_0), r(t_0)].
$$
We identify $\sigma (t_0, \nu ) \subset \Bbb C$ with its lift to
$W$ that contains $t_0.$
\medskip

{\bf Theorem A1.} {\it In the assumptions at the beginning of the
subsection, for any complex form $\omega \in \Omega_n^*,$}
\begin{equation}
\#\{t \in \sigma(t_0, \nu ) \mid I(t) = 0\}<(1-\log c')A^{578}, \ A
= e^{\frac{n^4}{c''}}(H). \label{1.5b}
\end{equation}
\medskip

This theorem is an immediate corollary of the Growth-and-Zeros
theorem and the Main Lemma stated below. Let
$$
L^{\pm}(t_0) = \begin{cases} \{a + \nu e^{\pm i \ph} \in W \mid \ph \in [0,
2\pi]\}
\text { for } a \not = -\infty\\
              \{-3 e^{\pm i \ph} \in W \mid \ph \in [0, 2(n+1)\pi]\},
\text { for } a = -\infty,\end{cases}
$$
$$
R^{\pm}(t_0) = \begin{cases} \{b - \nu e^{\pm i \ph} \in W \mid \ph \in [0,
2\pi]\}
\text { for } b \not = +\infty\\
              \{+3 e^{\pm i \ph} \in W \mid \ph \in [0, 2(n+1)\pi]\},
\text { for } b = +\infty,\end{cases}
$$
$$
\Gamma_a = L^+(t_0)\cup L^-(t_0),\ \Gamma_b=R^+(t_0)\cup R^-(t_0),
\
\Sigma =\Gamma_a\cup \Gamma_b\cup \sigma (t_0, \nu ). 
$$

\medskip
{\bf Main Lemma.}  {\it Let $H$    be a normalized
polynomial of degree $n+1\geq 3$ with critical values $
a_j: \ j= 1,..., n^2,$ $ \omega $ be a complex polynomial 1-form
of degree no greater than $ n.$ Let $W, \nu , \Sigma $ be the same
as at the beginning of this subsection. Then there exists a path
connected compact set $K \subset W$, \ $K\supset\Sigma$, $\pi
K\subset\overline{D_3}$, with the following properties:
\begin{equation}
diam_{int} K< 36n^2; \label{1.6}
\end{equation}
\begin{equation}
\rho(\pi K, a_j) \ge \nu \text { for any } j = 1,..., n^2.
\label{1.7}
\end{equation}
Moreover, let $U$ be the minimal simply connected domain in $W$
that contains the $\nu/2$ \nbd of $K.$ Then}
\begin{equation}
\bui< {(1 - \log c')}A^2. \label{1.8}
\end{equation}
\medskip


The Lemma is proved in Section 2. It is used also in the
estimate of the number of zeros of the integral in the intervals
$(a,l(t_0))$, $(r(t_0),b)$. In fact, a much better estimate for
the Bernstein index holds:
\begin{equation}
\bui < \frac{2700n^{18}}{c''(H)}-30n^6\log c'(H):= B(n, c', c'').
\label{bukpol}
\end{equation}
Inequality (\ref{bukpol})
is proved in 2.7. Together with the elementary inequality
\begin{equation}
B(n, c', c'') < (1 - \log c')A^2, \label{1.8a}
\end{equation}
it implies (\ref{1.8}).

\begin{proof} {\bf of Theorem A1.}
Let us apply Growth-and-Zeros theorem
to the function $I$ in the domain $U$ in order to estimate the
number of zeros of $I$ in $K;$ note that $K \supset \sigma (t_0).$
The intrinsic diameter of $K$ is estimated from above by (\ref{1.6}).
The gap condition for $U$ and $K$ has the form
$$
\e(K, \ \partial U) = \frac {\nu} 2 = \frac {c''} {8 n^2}
$$
by the definition of $U.$ Hence,
$$
e^{\frac{2D}{\e}}<e^{\frac{72n^2}{c''}8n^2} = A^{576}.
$$
The Bernstein index $\bui$ is estimated from above in (\ref{1.8}). By
Growth-and-Zeros theorem
$$
\#\{t \in \sigma(t_0) \mid I(t) = 0\} < B_{K,U}(I)A^{576} < (1 -
\log c')A^{578}.
$$

This proves (\ref{1.5b}).
\end{proof}

The following remark motivates the restriction on $R$ in Theorem
B.

\begin{remark} Let  $K$ be the set from the Main Lemma,
$\rho_WK$ be its diameter in the Poincar\'e metric of $W$. Then
\begin{equation}
\rho_WK<(c'')^{-1}288n^4.\label{1.9}
\end{equation}
Indeed, $\rho_WK$ is no greater than the ratio of the double
intrinsic diameter of $K$ divided by its minimal distance to the
critical values of $H$. Together with (\ref{1.6}) and (\ref{1.7}) this implies
(\ref{1.9}). On the other hand, in the proof of Theorem B, we apply
Growth-and-Zeros theorem in the case, when the Poincar\'e disc
$DP_R(H,t_0)$ is large enough, namely,  contains the set $K.$
\end{remark}

\subsection{Theorem A2 and proof of Theorem A}

{\bf Theorem A2.} {\it Let $H$, $t_0,a = a(t_0) , b = b(t_0)$ be the
same as in the previous subsection. Let $\omega$ be a real 1- form
in $\Omega_n^*$. Then, in assumptions of
Theorem A1,}
\begin{equation}
\#\{t \in (a,l(t_0))\cup(r(t_0),b)  \mid I(t) = 0\}< (1-\log
c')A^{4800} \label{1.9a}
\end{equation}
\medskip

\begin{proof} {\bf of Theorem A.}  By Theorems A1 and A2

\begin{equation}
\#\{ t\in(a,b),\ I(t)=0\} < (1-\log c')A^{578} +
(1-\log c')A^{4800} < 2(1-\log c')A^{4800}.  \label{1.10}
\end{equation}

This implies the estimate of the number of zeros given by Theorem
A on the interval $(a,b)$.

Let $\sigma'\subset\Bbb R$ be the maximal interval of continuity
of the family $\Gamma$ of real ovals that contains $\gamma_0.$
Then $\sigma'$ is bounded by a pair of critical values, at most
one of them may be infinite. In general, the interval $\sigma'$
may contain critical values (see Fig.1, which presents a possible
arrangement of level curves of $H$ in this case: $A_1$, $A_2$,
$A_3$ are critical points of $H$, $a_j=H(A_j)$, 
$a_2\in\sigma'=(a_1,a_3)$. 
In this case $\sigma'\neq(a,b)=(a_1,a_2)$. Let us estimate the
number of zeros on $\sigma'$. The interval $\sigma'$ is split into
at most $n^2$ subintervals bounded by critical values. On each
subinterval the number of zeros of $I$ is estimated by
(\ref{1.10}), as before. Therefore, the number of zeros of $I$ on
$\sigma'$ is less than $2n^2(1-\log c')A^{4800} < (1-\log
c')A^{4801}.$ This proves Theorem A.
\end{proof}

The paper is structured as follows. In Section 2 we prove the Main
Lemma modulo two statements: formula for the determinant of
periods, and upper estimates of Abelian integrals provided by
quantitative algebraic geometry. These two statements are treated
in two separate papers by the first author (A.Glutsyuk, \cite{G**}
and \cite{G*} respectively). After this the Main Lemma, as well as
Theorem A1, is proved. Theorem A2 is proved in Section 4. Theorem
B is proved in Sections 3 and 4. In both sections, the Main Lemma
is intensively used. The complete proof of Theorem A ends up in
Section 4.

\section{An upper bound for the number of zeros on a real
segment distant from critical values}

In this section we prove the Main Lemma and hence Theorem A1. We
also prove the Modified Main Lemma, see 2.9 below, and prepare
important tools for the proof of other results: Theorems A2, B1
and B2.

\subsection{Normalized ultra-Morse polynomials; notations}

Denote by $D_r$ a disk $|t| \le r.$

All  along this section $H$ is a real normalized ultra-Morse
polynomial of degree $n + 1 \ge 3, \ \mu = n^2; \ a_1, \dots ,
a_\mu $ are critical values of $H, \ \nu $ is the same as in
(\ref{1.2}), $\e = \nu/2.$ For $t$ close to $a_j, \ \delta_j(t)$ is a
local vanishing cycle corresponding to $a$ on a level curve
$$
S_t = \{ H = t\} ;
$$
the definition of this cycle is recalled in the next subsection.
Denote by $B = B_H$ the set of all noncritical values of $H:$
$$
B = \Bbb C \setminus \{ a_1. \dots , a_\mu \} .
$$
Let
$$
t_0 \in B\cap (-3, 3),
$$
and $W$ be the universal cover over $B$ with the base point $t_0$
and projection
$$
\pi : W \to B.
$$

\subsection{Marked system of vanishing cycles}

To begin, we recall well known results and definitions.

\begin{lemma} {\bf (Morse lemma).}
A holomorphic function having a Morse
critical point may be transformed to a sum of a nondegenerate
quadratic form and a constant term by an analytic change of
coordinates near this point.
\end{lemma}

\begin{corollary} Consider a holomorphic function in $\Bbb C^2$
having a Morse critical point with a critical value $a.$ An
intersection of a level curve of this function corresponding to a
value close to $a$ with an appropriate \nbd of the critical point
is diffeomorphic to an annulus.
\end{corollary}

This annulus may be called a local level curve corresponding to
the a critical value $a.$

\begin{definition} \label{d2.2} A generator of the first homology
group of the local level curve corresponding to $a$ is called a
{\it local vanishing cycle} corresponding to $a.$
\enddefinition
A local vanishing cycle is well defined up to change of
orientation.
\end{definition}
A path $\alpha_j: [0,1] \to \Bbb C$ is called {\it regular}
provided that
\begin{equation}
\alpha_j(0) = t_0, \ \alpha_j(1) = a_j, \ \alpha_j [0,1) \subset B
\label{2.1}
\end{equation}
\begin{definition} \label{d2.3} Let $\alpha_j$ be a regular path, $s \in
[0,1]$ be close to $1, \ \delta_j(t), \ t = \alpha_j(s),$ be a
local vanishing cycle on $S_t$ corresponding to $a_j.$ Consider
the extension of $\delta_j$ along the path $\a $ up to a
continuous family of cycles $\delta_j(s)$ in complex level curves
$H = \alpha_j(s).$ The homology class $\delta_j = \delta_j(0)$ is
called a {\it cycle vanishing along $\a_j .$}
\end{definition}
\begin{definition} \label{d2.4} Consider a set of regular paths $\alpha_1, \dots , \alpha_\mu
,$ see (\ref{2.1}). Suppose that these paths are not pairwise
and self
intersected. Then the set of cycles $\delta_j \in H_1(S_{t_0},\Bbb
Z)$ vanishing along $\alpha_j, \ j = 1. \dots , \mu ,$ is called a
{\it marked set of vanishing cycles} on the level curve $H = t_0.$
\end{definition}
\begin{definition} \label{d2.5} Any point $\hat t \in W$ is
represented by a class $[\lambda ]$ of curves in $B$ starting at
$t_0$ and terminating at $t = \pi \hat t;$ all the curves of the
class are homotopic on $B.$ Any cycle $\gamma $ from $H_1(S_{t_0},
\Bbb Z )$ may be continuously extended over $\lambda $ as an
element of the homology groups of level curves of $H;$ the
resulting cycle $\gamma (\hat t)$ from $H_1(S_t, \Bbb Z )$
is called an {\it extension} of $\g $ corresponding
to $\hat t.$
\end{definition}

Let  $\delta_1,..., \delta_\mu $ be a marked set of vanishing
cycles. For any  cycle $\delta_l$
from this set, denote by $W_l$ the Riemann surface of the integral
$$
I_l(\hat t) = \int_{\delta_l(\hat t)} \omega ,
$$
with the base point $t_0.$ Let $\pi_l$ be the natural projection
$W \to W_l.$ Denote by $D_r(a)$ the disk $|t - a| \le r.$

\begin{remark} The Riemann surface $W_l$ contains the disc
$D_{\nu}(a)$.
\end{remark}

\begin{lemma} {\bf (Modified Main Lemma).}
The Main Lemma from subsection 1.7 holds true provided that the
real oval $\gamma$ is replaced by a local vanishing cycle $
\delta_l(t)$ close to the corresponding critical value $a_l$,  and
$\Sigma $ is replaced by the disk $\overline D_{\nu}(a_l).$
\end{lemma}

This lemma is proved in 2.8.

\subsection{Matrix of periods}

Consider and fix an arbitrary marked set of vanishing cycles
$\delta_j, \ j = 1, \dots , \mu .$ For any $\hat t \in W,$ let
$\delta_j(\hat t)$ be the extension of $\delta_j$ corresponding to
$\hat t.$

\begin{definition} \label{d2.6} Consider a set $\Omega $ of $\mu $
forms $\omega_j$ of the type
\begin{equation}
\omega_i = yx^ky^ldx,\ k,l\geq0,\ k+l\leq2n-2
\label{2.2}
\end{equation}
$(k,l)$ depends on $i,$ such that all the forms with $k + l \le n
- 1$ are included in the set. In what follows, such a set is
called {\it standard}.
\end{definition}

A {\it matrix of periods} $\Bbb I = (I_{ij}), \ 1 \le i \le \mu ,
\ 1 \le j \le \mu $, is the matrix function defined on $W$ by the
formula:
\begin{equation}
I_{ij}(\hat t) = \int_{\delta_j(\hat t)}\omega_i, \ \Bbb I(\hat t)
= (I_{ij}(\hat t))   \label{2.3}
\end{equation}
where $\delta_j, \ j = 1, \dots , \mu,$  form a marked set of
vanishing cycles; $ \{ \omega_i| i = 1, \dots , \mu \} $ is a
standard set of forms (\ref{2.2}).

When we want to specify dependence on $H$, we write $\Bbb I(\hat
t,H )$ instead of $\Bbb I(\hat t).$

\subsection{Upper estimates of integrals}

Denote by $|\lambda |$ the length of a curve $\lambda ,$ and by
$U^\varepsilon (A)$ the $\varepsilon $-neighborhood of a set $A.$

The main result of the quantitative algebraic geometry that we
need is the following

\begin{theorem} \label{t2.1}
Let $\delta_j$ be a vanishing cycle from
a marked set, see Definition \ref{d2.4}, corresponding to a curve
$\alpha_j, |\alpha_j| \le 9$
(recall that $|t_0|\leq3$). Let $\lambda \subset B$ be a curve
starting at $t_0$ (denote $t$ its end) such that
\begin{equation}
|\lambda | \le 36n^2 + 1, \ |t|\leq4.   \label{2.4}
\end{equation}
Let the curve $\alpha_j\cap U^{\e}(a_j)$
be a connected arc of $\a_j$, and the curves
$\alpha_j \setminus U^\e (a_j)$ and $\lambda $
have an empty intersection with $\e $-\nbds of the critical values
$a_k,$ where $\e = \nu /2, \ \nu $ is from (\ref{1.2}).
Let $\omega $ be
a form (\ref{2.2}), $\hat t \in W$ corresponds to $[\lambda ],$ and
$\delta_j(\hat t)$ be the extension of $\delta_j$ to $\hat t$.
Then
\begin{equation}
|I_{\delta_j(\hat t)} \omega | < 2^{\frac{2600n^{16}}{c''(H)}}
(c'(H))^{-28n^4}:= M_0    \label{2.5}
\end{equation}
\end{theorem}
This result is based on Theorem $C$ from 1.5. Both results are
proved in the forthcoming paper \cite {G*}.

We have to give an upper bound of the integral not over a vanishing
cycle, but over a real oval. The following Lemma shows that
the real oval is always a linear combination of some
(at most $\mu$) vanishing
cycles with coefficients $\pm1$.

\def\gto{\gamma_{t_0}}

\begin{lemma} \label{l1.3} {\bf (Geometric lemma).}
Let $H$ be a real ultra-Morse polynomial and $\g $ be a real oval
of $H.$ Let $H|_\g = t_0.$ Consider the critical values of $H$
that correspond to the critical points located inside $\g $ in the
real plane. Let $\alpha_j, \ j = 1, \dots , s$, be nonintersecting 
and nonself-intersecting paths that connect $t_0$ with these
critical values and satisfy assumption (\ref{2.1}); we may change
the numeration of critical points to get the first $s$ ones inside
$\g.$ Moreover, suppose that all these paths belong to the upper
halfplane and no open domain bounded by a path $\alpha_j$ and a
real segment (connecting the endpoints of $\alpha_j$) contains any
critical value of $H$ (see Fig.2a,b).  
Let $\delta_j$ be the vanishing cycles that correspond to the paths
$\alpha_j.$ Then
\begin{equation}
[\g ] = \Sigma_1^s\e_j\delta_j, \text { where } \e_j =
\pm 1.\label{1.5}
\end{equation}
\end{lemma}

\begin{figure}[ht]
  \begin{center}
\label{fig:2}
   \epsfig{file=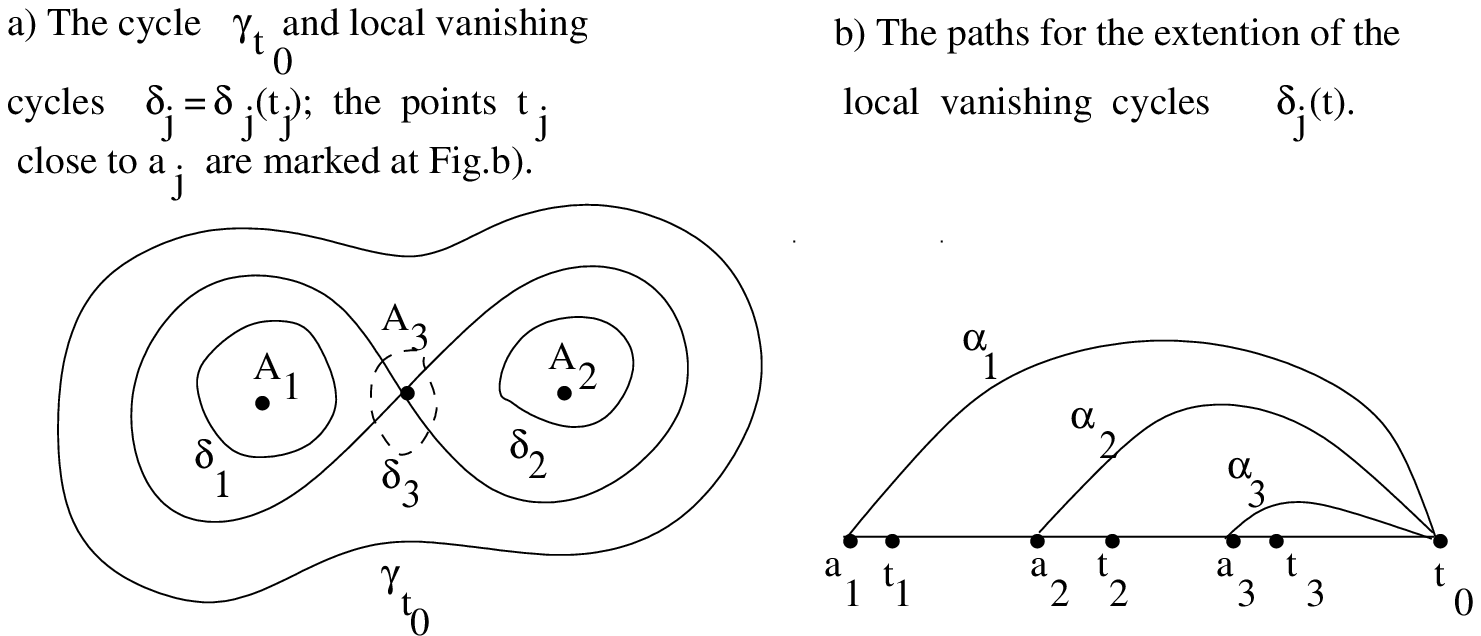}
    \caption{}
  \end{center}
\end{figure}

The authors believe that Lemma \ref{l1.3} is well known to
specialists, but they did not find it in literature. Its proof
is given in 3.5.

\begin{corollary} \label{cor1} The integral (\ref{1.1}) constructed for the
real oval $\g = \gamma (t_0)$ satisfies the upper estimate:
$$
\left| \int_{\g (\hat t)}\omega \right| \le n^2\max_{j=1,\dots
,s}\left| \int_{\delta_j(\hat t)}\omega \right| .
$$
\end{corollary}

\begin{corollary} In the condition of the previous Theorem let
$H$ be a real polynomial, $\gamma(\hat t)$ be the extension to
$\hat t$ of a real oval,
\begin{equation}
\omega=\sum_{k+l\leq n-1}a_{kl}yx^ky^ldx.\label{omkl}\end{equation}
Then
\begin{equation}|I_{\gamma(\hat t)}\omega|\leq n^4M_0
\max_{k+l\leq n-1}|a_{kl}|.\label{2.5'}
\end{equation}
\end{corollary}
\subsection{Determinant of periods}

The determinant of the matrix of periods (\ref{2.3})
is called the {\it
determinant of periods.} It appears that this determinant is
single-valued on $B,$ thus depending not on a point of the universal
cover $W,$ but rather on the projection of this point to $B.$ Let
$$
\Delta (t) = \det \Bbb I(\hat t), \ t = \pi \hat t.
$$
The single-valuedness of the main determinant follows from the
Picard-Lefschetz theorem. Indeed, a circuit around one critical
value adds the multiple of the correspondent column to some other
columns of the matrix of periods. Thus the determinant remains
unchanged.

When we want to specify the dependence of the main determinant on
$H,$ we write $\Delta (t).$ This function is polynomial in $t,$
and an algebraic function in the coefficients of $H.$ The formula
for the main determinant (with $\omega_i$ of appropriate degrees)
with a sketch of the proof was claimed by
A.Varchenko \cite {V89}; this formula is given up to a constant
factor not precisely determined. The complete answer
(under the same assumption on the degrees of $\omega_i$)
is obtained by the first author (A.Glutsyuk, \cite {G**}).
Moreover, the following lower estimate holds:

\begin{theorem} \label{t2.2}
For any normalized ultra-Morse polynomial
$H,$ the tuple $\Omega $ of standard forms (\ref{2.2}) may be so chosen
that for any $t\in\mathbb C$ lying outside the $\nu=
\frac{c''}{4n^2}$- neighborhoods of the critical values of $H$
the following lower estimate holds:
\begin{equation}
|\Delta (t,H)| \ge {(c'(H))}^{6n^3}{(c''(H))}^{n^2}n^{-62n^3}:=
\Delta_0   \label{2.6}
\end{equation}
\end{theorem}

This result is proved in \cite {G*} with the use of the explicit
formula for the Main Determinant mentioned before, and results of
the quantitative algebraic geometry.

\subsection{Construction of the set $K$}

We can now pass to the construction of the set $K$ mentioned in
the Main Lemma. We first construct a smaller set $K'.$

\begin{lemma} \label{l2.1} {\bf (Construction lemma).}
Let $\gamma \subset S_{t_0}$ be a real oval.
There exists a set of regular paths
$\alpha_j, \ j = 1, \dots , \mu ,$ see Definition \ref{d2.3},
such that:
$$
|\alpha_j| \le 9,
$$
the paths $\alpha_j$ are not pairwise and self intersected;

and there exists a path connected set $K' \subset W, \ t_0 \in K',
\ \pi K' \subset D_3,$ such that:

for any cycle $\delta_j \in H_1(S_{t_0}, \Bbb Z)$ vanishing along
$\alpha_j$ there exist two points $\tau_1, \tau_2 \in K' \cap \pi
^{-1}(t_0)$ such that
\begin{equation}
[\gamma (\tau_1)] - [\gamma (\tau_2)] = l_j[\delta_j],\
l_j\in\mathbb Z\setminus0.   \label{2.7}
\end{equation}
Moreover,
\begin{equation}
{\text {diam}}_{\text {int}}K'<19n^2,   \label{2.8}
\end{equation}
and $\pi K'$ is disjoint from $\nu $-\nbds of the critical values
$a_j, \ j = 1, \dots , \mu .$
\end{lemma}

\begin{lemma} \label{l2.2} {\bf (Construction lemma
for vanishing cycles).}
Construction lemma holds true if $\g  \subset S_{t_0}$ is replaced
in its statement by any vanishing cycle $\delta_l = \delta_l(t_0)$
from an arbitrary marked set of vanishing cycles, and $W$ is
replaced by $W_l.$ In the conclusion, (\ref{2.7})
should be replaced by
$$
[\delta_l (\tau_1)] - [\delta_l(\tau_2)] = l_j[\delta_j], \text { for
} j \not = l, \ \ [\delta_l (\hat t)] = [\delta_l] \text { for
}\hat t = t_0,\ l_j\in\mathbb Z\setminus0.
$$
\end{lemma}

Both lemmas are  proved in 2.9.  In what follows we deduce the
Main Lemma from Lemma \ref{l2.1} and Theorems \ref{t2.1}, \ref{t2.2}.

\begin{corollary} {\bf (of Lemma \ref{l2.1}).}
For any form $\omega $ (not necessarily of type
(\ref{2.2})) and any marked set of vanishing cycles consider the
vector function
\begin{equation}
\mathbb I_\omega : W \to {\Bbb C}^\mu , \ \hat t \mapsto \left
(\int_{\delta_1(\hat t)}\omega , \dots , \int_{\delta_{\mu} (\hat
t)}\omega \right ).   \label{2.9}
\end{equation}
Let $||\cdot ||$ denote the Euclidean length in $\Bbb C^\mu .$
Then
\begin{equation}
m_0:= \max_{\hat t \in K'\cap \pi ^{-1}(t_0)}|I(\hat t)| \ge \frac
{1}{2n}||\mathbb I_\omega (t_0)||.  \label{2.10}
\end{equation}
\end{corollary}

\begin{proof} Let us take $j$ so that
$$
\max_i\left |\int_{\delta_i(t_0)}\right | = \left
|\int_{\delta_j(t_0)}\right |.
$$
Then
\begin{equation}
\left |\int_{\delta_j(t_0)}\right | \ge \frac {1}{n}||I_\omega
(t_0)|| \label{2.11}
\end{equation}
By Lemma \ref{l2.1}, there exist $\tau_1, \tau_2$ such that
$$
I(\tau_1) - I(\tau_2) = l\int_{\delta_j(t_0)}\omega , \ l \in \Bbb
Z \setminus 0.
$$
Hence, at least one of the integrals in the left hand side, say,
$I(\tau_l), \ l \in \{ 1,2\} ,$ admits a lower estimate:
\begin{equation}
|I(\tau_l)| \ge \frac {1}{2}\left |\int_{\delta_j(t_0)}\omega
\right |. \label{2.12}
\end{equation}
Together with (\ref{2.11}) this proves the corollary.
\end{proof}

Let us now take
\begin{equation}
K = K' \cup \Sigma , \Sigma = \sigma (t_0)\cup L^\pm (t_0)\cup
R^\pm (t_0).  \label{2.13}
\end{equation}

In the following section we will check that this $K$ satisfies the
requirements of the Main Lemma.

\subsection{Proof of the Main Lemma}

Let us take $K$ as in (\ref{2.13}).
Let $\nu $ be the same as in (\ref{1.2}).
Let $U$ be the smallest simply connected set that contains the $\e
$-\nbd of $K, \ \e = \nu /2.$ Then (\ref{1.6}) follows from
(\ref{2.8}),
(\ref{2.13}). The last statement of Lemma \ref{l2.1}
implies (\ref{1.7}).

Let us now check (\ref{1.8}), that is, estimate
from above the Bernstein
index $B_{K,U}(I).$

By Theorem \ref{t2.1} and (\ref{1.6}),
all the elements of the matrix $\Bbb I(t_0)$ admit an
upper bound:
$$
|I_{ij}(t_0)| < M_0.
$$
Fix a form $\omega_0 = A_ndx + B_ndy.$ There exists another form
$\omega$ of the type (\ref{omkl}) such that the
form $\omega - \omega_0$ is exact. Let $\omega_i =
yx^{k_0}y^{l_0}dx$ be such that $|a_{k_0l_0}| = \max_{k+l\le
n-1}|a_{kl}|$ in (\ref{omkl}).
Without loss of generality we set $a_{k_0l_0} = 1.$
Let us now replace the $i$th row of the matrix $\Bbb I$ by the
vector $\Bbb I_\omega .$ This transformation is equivalent to
adding to the $i$-th line linear combination of other lines, so the
determinant $\Delta (t_0)$ remains unchanged. All the elements in
all other rows are estimated from above by $M_0.$ Hence, all the
vector-rows except for the $i$th one have the length at most
$nM_0.$ By (\ref{2.10}), the $i$th  row has the length at most 
$2nm_0.$
Hence,
$$
\Delta_0 \le 2m_0M_0^{\mu -1}n^\mu ,\ \mu=n^2
$$
where $\Delta_0$ and $M_0$ are from (\ref{2.6}) and (\ref{2.5})
respectively.
Therefore,
\begin{equation}
\log m \ge \log m_0\ge 
\log \Delta_0 - (\mu - 1) \log M_0 - \mu \log n-\log2.
\label{2.12a}\end{equation}

On the other hand, by (\ref{2.5'}),
$$
\max_U|I| \le n^4M_0.
$$
Hence,
$$
B_{K,U}(I) \le  (\mu + 4)\log n + \mu \log M_0 - \log \Delta_0+\log2.
$$
Now, elementary estimates imply (\ref{bukpol}).
This proves the Main Lemma.

\subsection{Modified Main Lemma and zeros of integrals over
(complex) vanishing cycles}

\begin{proof} {\bf of the Modified Main Lemma.} The arguments of the
previous section work almost verbatim. The previous Corollary for
$\gamma$ replaced by $\delta_l$ is stated and proved in the same
way.

Let $K'$ be the same as in Lemma \ref{l2.2}. Instead of
(\ref{2.13}), let
$$
K = K' \cup \a_l \cup\overline D_{\nu} (a_l).
$$
Let $U$ be the smallest simply connected set that contains the
$\e$-\nbd of $K.$

By Theorem \ref{t2.1},
$$
\max_{\overline V} |I_l | \le M_0, \text { where } V = U \setminus D_\nu
(a_l).
$$

But $I_l$ is holomorphic in $ D_\nu (a_l).$ Hence, by the maximum
modulus principle, the previous inequality holds in $U$ instead of
$V.$ After that, the rest of the arguments of the previous section
work. This proves the Modified Main Lemma.
\end{proof}

\begin{theorem} \label{t2.3} The number of zeros of the
integral $I_l$ in the disk $ D_\nu (a_l)$ satisfies the inequality:
\begin{equation}
\#\{ \hat t \in D_\nu (a_l)| I_l(\hat t) = 0 \} \le (1 - \log
c'(H))A^{578}.       \label{2.14}
\end{equation}
\end{theorem}

The proof is the same as for Theorem A1,  section 1.7.


\subsection{Proof of the Construction Lemmas}


\begin{proof} {\bf of Lemma \ref{l2.1}.}

\begin{definition} \label{d2.7} A loop $\lambda_j$  is associated to a regular path
$\a_j$ if
$$
\lambda_j=\tilde\a_j'\partial D_{\nu}(a_j)(\tilde\a_j')^{-1},
$$
where $\tilde\a_j'=\a_j\setminus D_{\nu}(a_j), \ \nu=\frac{c''}{4n^2},\
\partial D_{\nu}(a_j)$ is positively oriented
(we suppose that $\tilde\a_j'$ is connected).
\end{definition}

Let $\alpha_1, ..., \alpha_{\mu} $ be the same as in Definition
\ref{d2.4}.
The set $K'$ we are looking for will be the union of appropriate
$n^2$ liftings of the loops $\lambda_j$ (one lifting for each
$\lambda_j$) associated with $\alpha_j$ to the Riemann surface $W.$
In what follows, the choice of the curves $\alpha_j$ will be specified.

We prove Lemma \ref{l2.1} in four steps. The set $K'$
is constructed in the first three steps. In the fourth step
we check that the resulting set has the
required properties.

Step 1: special path set. Denote by $\a'_j$ the segment
$[t_0,a_j]$ oriented from $t_0$. Fix $j$ and suppose that $\a_j'$
contains critical values of $H$ different from $a_j$; denote the
set of these values by $A$. For any $a_i\in A$ replace the
diameter $\a_j'\cap D_{\nu}(a_i)$ by a semicircle.

If $a_j$ is real, then this semicircle is chosen in the upper
half-plane. In general, for any fixed line $\theta$ passing
through $t_0$ and some critical values the previous semicircles
corresponding to all the pairs $a_i,a_j$ in $\theta$ as above are
chosen to be on one and the same side from $\theta$. If $\a_j'$
intersects a disc $D_{\nu}(a_s)$ but does not contain $a_s$,
replace the chord $\a_j'\cap D_{\nu}(a_s)$ by the smallest
arc of the circle $\partial D_{\nu}(a_s). $
The path thus constructed will be denoted by $\a_j$. Recall
that $t_0\in\bar{D_3}$, $a_j\in D_2$. Therefore, the length
of any segment $\a_j'$ is less than 5. Hence,
\begin{equation}
|\a_j|\leq\frac{5\pi}2<9. \label{2.15}
\end{equation}

Each path $\a_j$ is nonself-intersected by construction and is
contained in $D_3$ (except may be for $t_0$). One can achieve that
the paths $\a_j$ be disjoint outside $t_0$ by applying to them
arbitrarily small deformation preserving the previous inequality
and inclusion.

Step 2: special loop set. For any $j$ denote  by $\lambda_j$ the
loop associated to $\a_j$ in the sense of Definition \ref{d2.7}.
By construction, $\lambda_j\subset\bar{D_3}$. We have
$$
|\lambda_j|\leq2|\a_j|+|\partial D_{\nu}(a_j)|<19.
$$

\def\gi {\delta_i}
\def\gj {\delta_j}

Step 3: construction of $K'$. Denote by $G$ the intersection graph 
of $\gto$ and all the vanishing cycles $\gi$ (along the previously
constructed paths $\a_i$). {\it This graph is connected.}
This follows from the two lemmas below.

\begin{lemma} \label{l2.3} The intersection graph of the marked set of vanishing
cycles is connected. The set itself forms a basis in the group $H_1(S_{t_0}, \Bbb Z).$
\end{lemma}

(Recall the definition
of the intersection graph: its vertices are identified with the cycles;
two of them are connected by an arc, if and only if the corresponding
intersection index is nonzero.)

Lemma \ref{l2.3} is implied by  the following statements from
\cite{1}: theorem 1
in 2.1 and theorem 3 in 3.2.

\begin{lemma} \label{l2.4} Consider a maximal family of real ovals that
contains $\gto.$ The union of the ovals of the family forms an
open domain. The boundary of this domain consists of one or two
connected components. Any of these components belongs to a
critical level of $H$ and contains a unique critical point. Fix
any of these critical points and denote by $\delta $ the
corresponding local vanishing cycle. Then the cycle $\de$ may be
extended to a cycle $\de (t_0)$ that belongs to a marked set of
vanishing cycles constructed above. Moreover,
$$
(\de (t_0), \gto) \not = 0, \ \text{more precisely, it is equal to}
\ \pm1, \ \pm2.
$$

\end{lemma}

The proof of this lemma is written between the lines of \cite{9}, pp 12,13.
It is illustrated by Fig.3.

\begin{figure}[ht]
  \begin{center}
   \epsfig{file=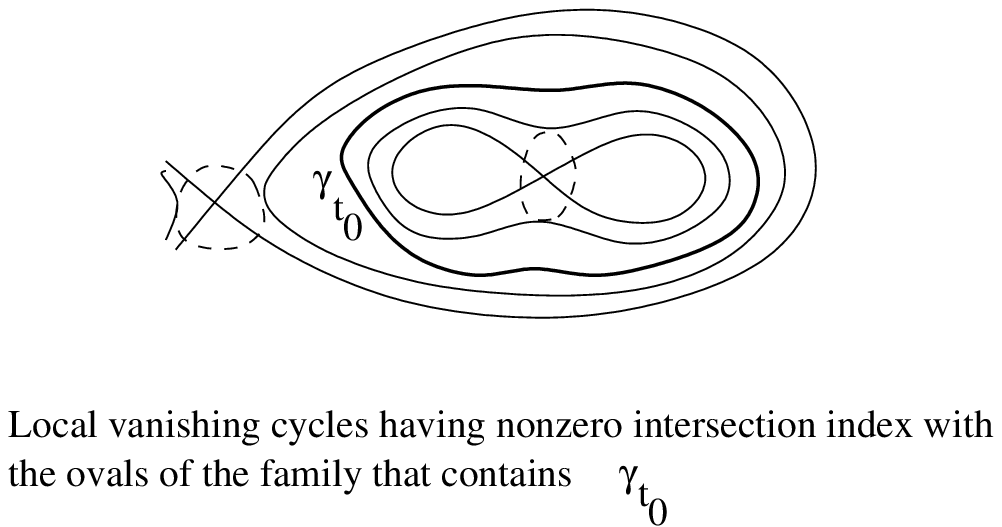}
    \caption{}
    \label{fig:3}
  \end{center}
\end{figure}

Let us define a metric on the set of the vertices of the graph $G$.
Suppose that each edge of $G$ has length 1. Then the distance
$D_G$ between any two vertices of $G$ is well defined as the
length of the shortest path in $G$ that connects the vertices. For
any $r\in \Bbb N$ let
$$
S_r=\{\delta_j\ | \ D_G(\gto,\delta_j)=r\}.
$$
Let $T$ be a maximal tree in $G$ with the root $[\gto]$ such that
the distance in $T$ (defined as $D_G$ but with paths in $T$) of
any vertex to the root $[\gto]$ coincides with $D_G$ (see Fig.4, 
where the tree $T$ is marked by bold curves.)

\begin{figure}[ht]
  \begin{center}
   \epsfig{file=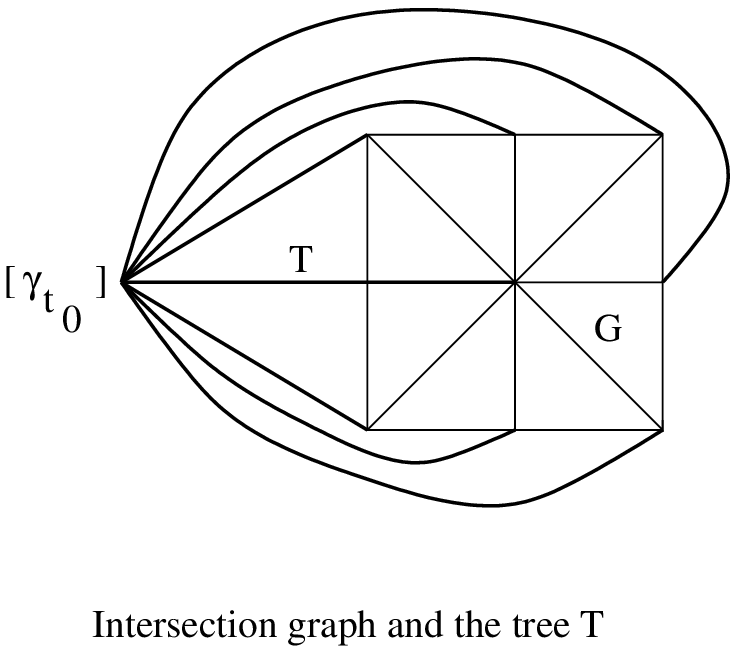}
    \caption{}
    \label{fig:4} 
  \end{center}
\end{figure}

For any vanishing cycle $\delta_j(t_0)$ let $L_j$ be the branch of
the tree $T$ from $[\gto]$ to $\gj(t_0)$. Let $[\gto]$,
$\delta_{j_1}(t_0)$, \dots, $\delta_{j_r}(t_0)=\delta_j(t_0)$ be
its vertices ordered from the beginning to the end of the branch.
By definition, the intersection index of any cycle in this
sequence with its two neighbors is nonzero, and that of any two
nonneighbor cycles is zero. Let us call this {\it the regularity
property} of $L_j$.

The set $K'\subset W$ we are looking for is the image of the tree
$T$ in $W$ under a continuous map $\phi:T\to W$. This map is
defined by induction in $r$ as follows. It suffices to define
$\phi|_{L_j}$ for any $\delta_j$.

Base of induction: $r=0$. The cycle $\gto$ is mapped  to $t_0$.

Induction step. Suppose that the cycle $\delta'=\delta_{j_{r-1}}$
is mapped to $\tau_1\in\pi^{-1}(t_0)$: $\phi(\delta')=\tau_1$. Let
us lift the loop $\lambda=\lambda_{j_r}$ to $W$ as a covering
curve $\tilde{\lambda}$ over $\lambda$ with the starting point
$\tau_1$. Let $\delta=\delta_{j_r}$, $\tau_2\in\pi^{-1}(t_0)$ be
the endpoint of $\tilde{\lambda}$. This induces a map of the edge
$[\delta',\delta]$ to $\tilde{\lambda}$. This map defines the
extension of $\phi$ to the edge $[\delta',\delta]$. The induction
step is over.

Step 4: properties of the set $K'$. The set $K'$ is a curvilinear
tree and thus, path connected. Its intrinsic diameter admits the
upper estimate
$$
{\text {diam}}_{int}K'\leq n^2\max_j|\lambda_j|<19n^2.
$$
The set $K'$ is projected to the loops $\lambda_j$, which lie in
$\bar{D_3}$ and are disjoint from the $\nu$- neighborhoods of the critical values by
definition. Hence, the same is true  for $\pi (K')$.

For any cycle $\delta=\delta_j(t_0)$ vanishing along the path $\a_j$ from
the special path set, see Step 1, let $L$ be the edge of the tree $T$ with
the endpoint $\delta$. Let $\delta'$ be the initial point of $L$. Let
$\tau_1=\phi(\delta'),\ \tau_2=\phi(\delta)$. Then (\ref{2.7}) holds   by the
Picard-Lefschetz theorem. In more details, let $L_j$ and $\delta_{j_m}(t_0)$
be the same, as in Step 3. Then
\begin{equation}
\gamma_{\tau_2}=\gto+\sum_{m=1}^rl_m\delta_{j_m}(t_0),\
l_m\in\Bbb Z\setminus 0, \label{2.16}
\end{equation}
\begin{equation}
\gamma_{\tau_1}=\gto+\sum_{m=1}^{r-1}l_m\delta_{j_m}(t_0),\
l_m\in\Bbb Z\setminus0. \label{2.17}
\end{equation}

Let us prove (\ref{2.16}) by induction in $r$  taking
(\ref{2.17}) as the induction hypothesis. Equality (\ref{2.17})
implies (\ref{2.16}) by Picard-Lefschetz theorem \cite{1} and the
regularity property of $L_j$, see Step 3. On the other hand,
(\ref{2.16}) and (\ref{2.17}) imply (\ref{2.7}). Lemma
\ref{l2.1} is proved.
\end{proof}

Lemma \ref{l2.2} is proved in the same way with the following
minor changes:
$G$ is now the intersection graph of the marked
set of vanishing cycles concidered, and in the lifting
process, $W$ should be replaced by $W_l.$

\section{Number of zeros of abelian integrals in complex
domains distant from critical values}

In this section we prove the first part of Theorem B, namely,

\subsection{Upper estimates in Euclidean and Poincar\' e disks}

Theorem B1 mentioned in 1.5 is stated below.

All through this section notations of section 2.1 hold. Moreover,
$\delta_1, \dots , \delta_{\mu} $ is a marked set of vanishing
cycles on $S_{t_0}, \ K \subset W$ is a compact set from the Main
Lemma, see 1.7.

\begin{theorem} \label{t3.1} 
Fix a normalized polynomial $H.$ Let $\hat t \in W$ be a point
represented by a curve $\lambda \subset B.$ Let $\delta_j$ be a
vanishing cycle from a marked set corresponding to a curve
$\alpha_j.$ Let $\alpha =\lambda^{-1}\alpha_j.$ Let 
$0<\beta<1$. Suppose that $\a\cap D_{\beta}(a_j)$ is a connected 
arc of the path $\a$, and $\a$ avoids the $\beta$- neighborhoods of 
the critical values distinct from $a_j$ of the polynomial $H$. 
Then for any 1- form $\omega $ of type (\ref{2.2}) 
\begin{equation}
\left| \int_{\delta_j(\hat t)}\omega \right|< 
2^{-2n} M_1,\ M_1=2^{10n^{12}\frac{|\a|+5}{\beta}}(c'(H))^{-28n^4} 
\label{3.1}
\end{equation}
\end{theorem}
Theorem \ref{t3.1} is proved in \cite {G*}. It is used in the 
estimate of the number of zeros  in Euclidean disc. The 
following upper bound (Theorem \ref{tpoinc} also proved in 
\cite{G*})  
of integrals is used to prove an upper bound of the number 
of zeros in Poincar\'e disc that is exponential in the 
radius of the disc.

\begin{remark} One can estimate the number of zeros in 
Poincar\'e disc by using Theorem \ref{t3.1} instead of 
Theorem \ref{tpoinc} (see the proof for Euclidean disc below). 
But the upper bound of the number 
of zeros obtained in this way is double exponential in the 
radius. 
\end{remark}

\def\var{\varepsilon}
\def\la{\lambda}
\def\zz{\mathbb Z}
\def\Var{\operatorname{V}}
\def\wt#1{\widetilde#1}
\def\cc{\mathbb C}

Denote by $\Var_{\gamma}f$ the variation of the argument of the function 
f along an oriented curve $\gamma$. 

\begin{theorem} \label{tpoinc} \cite{G*}. 
Let $H$ be a normalized ultra-Morse 
polynomial of degree $n+1\geq3$. Let $\hat t\subset W$ be a 
point represented by a curve $\lambda\subset B$. Let 
$\delta$ be a vanishing cycle 
from the marked setp; $\delta$ corresponds to a curve $\a_j$, 
$\a=\la^{-1}\a_j:[0,1]\to B$. 
Let $0<\beta\leq\nu=\frac{c''(H)}{4n^2}$, 
$$t'=\a(0)=\pi(\hat t),\ a=\a(1)=\a_j(1),\ \tau'=\min\{\tau\in[0,1]\ | \ \a(\tau,1]
\subset D_{\beta}(a)\}, \ \hat\a=\a\setminus\a(\tau',1],$$
\begin{equation}\wt\a=\a\cap 
(\overline D_3\setminus\cup_iD_{\beta}(a_i)),\label{cald}\end{equation}
$$
V=V_{\a,\beta}=\beta\sum_i\Var_{\hat\a\cap D_{\beta}(a_i)}(t-a_i)
+3\Var_{\hat\a\setminus\overline D_3}t.$$
Let $\delta\in H_1(S_{t'},\zz)$ be the cycle vanishing along $\a$. 
 Let $\omega$ be a monomial 
1- form of degree at most $2n-1$ with unit coefficient. Then 
\begin{equation} |I_{\delta}(\hat t)|\leq 2^{-2n}M_2,\ M_2=
2^{20n^{12}\frac{|\wt\a|+V+5}{\beta}}
(c'(H))^{-28n^4}\max\{1,(\frac{|t'|}5)^2\}.\label{uppi3}\end{equation}
\end{theorem}

Let $DE_{R,\beta }$ be an Euclidean disk in $W$ with $\beta $-
neighborhoods  of critical values deleted. More precisely,
$DE_{R,\beta }$ is the set of all those $\hat t \in W$ that may be
represented by a curve $\lambda ,$ whose length is no greater than
$R,$ provided that $\lambda $ avoids $\beta $-\nbds of critical
values.

\begin{theorem} \label{t3.2}  
Let $H$ be a normalized complex ultra-Morse 
polynomial of degree $n+1\geq3$, $\omega$ be arbitrary 1- form
of degree at most $n$. Then the number of
zeros of integral (\ref{1.1}), which is an analytic extension of an
integral over real ovals or over marked vanishing cycles of a
normalized polynomial $H$, is estimated from above as follows:
\begin{equation}
\# \{ \hat t \in DE_R|I(t) = 0\} <  
(1 - \log c'(H))e^{\frac{9R}{\beta}}   \label{3.2}
\end{equation}
provided that
\begin{equation}
R\ge 36n^2, \ \beta \le \nu/2    \label{3.3}
\end{equation}
\end{theorem}

The following statement is an analogue of Theorem \ref{t3.2}
for Euclidean metric replaced by the Poincar\' e one.

\begin{theorem} \label{tb1} {\bf (Theorem B1).}
In the assumptions of Theorem \ref{t3.2} the number of
zeros of integral (\ref{1.1}) over real ovals or over marked 
vanishing
cycles of a normalized polynomial $H$ is estimated as follows:
\begin{equation}
\# \{ \hat t \in DP_R|I(\hat t) = 0\} <(1-\log c'(H))e^{7R},
\label{3.4}
\end{equation}
provided that
\begin{equation}
R \ge \frac {288n^4}{c''(H)}.   \label{3.5}
\end{equation}
\end{theorem}
Recall that $DP_R$ is the disk in the Poincar\' e metric of $W$
of radius $R$ centered at the base point $t_0.$

Theorem B1 forms the first part of Theorem B. The second part of
Theorem B, Theorem B2, is presented in Section 4. Theorems B1, B2
imply Theorem B.

\subsection{Idea of the proof}

Theorems \ref{t3.2} and B1 are proved as Theorem A1, making use of
Growth-and-Zeros Theorem. The set $K,$ both from the Main Lemma
and from the Modified Main Lemma, belongs to $DE_R$ by (\ref{3.3})
and to $DP_R$ by (\ref{3.5}), see (1.11) and Remark 1.15.

Thus we have the main ingredient in the estimate of the Bernstein
index, namely, the lower bound for $m,$ see (\ref{2.12a}).

An upper estimate for the integral over a vanishing cycle is
provided by Theorem \ref{t3.1}. Yet there is a gap to
be filled when we
wish to replace a vanishing cycle by a real oval. This is done
by the following corollary of the Geometric Lemma \ref{l1.3}.

\begin{corollary} \label{cor2}
In the assumptions of Theorem \ref{t3.1} and any real
oval $\g $ of $H,$ 
\begin{equation}
\left| \int_\g \omega \right| < \frac{M_1}{n^2},   \label{3.6}
\end{equation}
see (\ref{3.1}).
\end{corollary}

Now everything is ready for the application of the Growth-and
Zeros theorem.

\subsection{Number of zeros in a Euclidean disk}

\begin{proof} {\bf of Theorem \ref{t3.2}.}
Denote the closure of the domain
$DE_R$ by $K.$ Let $\e' = \beta/2$,
$U$ be the
smallest simply connected domain in $W$ that contains the $\e'
$-\nbd of $K.$ Then
$$
D:= {\text {diam}}_{int}K \le 2R, \ \pi \text {-gap }(K,\partial
U) = \e' .
$$
Hence,
$$
e^{\frac{2D}{\e }} \le e^{\frac{8R}{\beta}}.
$$
This is the main factor in the estimate (\ref{3.2}).

Let us now estimate from above the Bernstein index $B =
B_{K,U}(I).$ Let $K'$ be the set from Lemma \ref{l2.1} (case of
real oval) or Lemma \ref{l2.2} (case of vanishing cycle). Let
$$
m = \max_K|I|.
$$
One has $K'\subset K$ by (\ref{2.8}) and
(\ref{3.3}). Therefore,
$\log m \ge \log m_0,$ where $\log m_0$ is from
(\ref{2.10}). On the other hand, let
$$
M = \max_U|I|.
$$
As in the proof of the Main Lemma, we assume
(without loss of generality) that $\omega$
is of the type (\ref{omkl}) with $\max|a_{kl}|=1$. Then
by Corollary \ref{cor2} (case of real cycle) or
Theorem \ref{t3.1} (case of
vanishing cycle), one has
$$
M \le M_1,
$$
$M_1$ is from (\ref{3.1}). Then
$$
B_{K,U}(I) \le \log M_1 - \log m_0.
$$
Inequalities (\ref{2.12a}), (\ref{3.3}) together with 
elementary estimates imply that
$$
\log M_1 - \log m_0 < 
(1 - \log c'(H))e^{\frac R{\beta}}.
$$
Together with
Growth-and-Zeros theorem, this completes the proof of Theorem
\ref{t3.2}.
\end{proof}

\subsection{Number of zeros in a disk in the Poincar\' e metric}

The proof of Theorem B1 is carried on by application of 
version (\ref{1.5n}) of the Growth-and-Zeros Theorem to the sets 
$$K_R=\overline{DP_R},\ \ U_R=DP_{R+1}:$$
\begin{equation}\#\{\hat t\in\overline{DP_R}\ | \ I(\hat t)=0\}\leq B_{K_R,U_R}
e^{\rho_R},\ \rho_R=diam_{PU_R}K_R.\label{krur}\end{equation} The
right-hand side of the latter inequality is estimated below. 

The set $K$ from the Main Lemma is contained in $K_R$ 
(this follows from (\ref{1.9})), and as 
before, this yields immediately lower bound of $m$. 
The principal part of the proof of Theorem B1 is the upper
bound of the integral on the set $U_R$:
\begin{equation}
max_{\overline U_R}|I(\hat t)|<M(R),\ \log M(R)=
(1-\log c'(H))e^{1.2R}.
\label{maxir}\end{equation}
 To prove it, we use
Theorem \ref{tpoinc}. Namely, given a $\hat t\in U_R$, consider
the path $\hat\lambda$ that is the geodesic from $t_0$ to $\hat t$
in the Poincar\'e metric of $W$ (we put $\la=\pi(\hat\la)$) 
and the path 
$$\a=\la^{-1}\a_i \ \text{from}\ t=\pi(\hat t) \ 
\text{to}\ a_i.\ \text{Put}\ \beta=\nu=\frac{c''}{4n^2}.$$ 
We have to estimate from above the value $M_2$ 
from Theorem \ref{tpoinc}, in particular, to estimate from above 
the module $|t|$ and the linear combination $V_{\a,\nu}$ of 
variations.  
To do this, we prove the following upper bound of the radius of 
the closed Euclidean disc containing 
$\overline{DP_{R+1}}\supset\la$ and lower bound of the gap between 
$DP_{R+1}$ and the critical values of $H$:

\begin{equation} \pi(\overline{DP_{R+1}})\subset D_{M_R}, \ \text { where } 
\log M_R=6e^R\log R,\
\label{mr}\end{equation}

\begin{equation}
dist(\overline{\pi(DP_{R+1}}),a_i)>\beta_R,\ \text { where } \beta_R=
M_R^{-1}.\label{betar}\end{equation}

Using the two latter inequalities we show that 

\begin{equation} V_{\a,\nu}<37e^RR\log R.\label{var<<}\end{equation}

The proofs of (\ref{mr})-(\ref{var<<}) 
and upper bound of $|\wt\a|$ ($\wt\a$ is defined in (\ref{cald}))  
are based on the 
following lower bounds of the Poincar\'e metric. Given a domain 
$G\subset\cc$, $\#(\cc\setminus G)>1$,  
denote by $P(G)$ the ratio of the Poincar\' e metric of $G$ to the
Euclidean one; $P(G)$ is a function in $t \in G.$

\def\cc{\mathbb C}

{\bf Inequality} (follows from theorem 2.17 in \cite{32}). For any
distinct $a,b\in\mathbb C$ one has
\begin{equation} P(\cc\setminus\{ a,b\})(t) >
[\min_{c=a,b}|t-c|(\min_{c=a,b}|\log|\frac{t-c}{a-b}||+5)]^{-1}
\label{pab}\end{equation}

\begin{corollary} Let $H$ be a balanced polynomial, $B$ be
the complement of $\cc$ to its critical values. Then
\begin{equation} P(B)(t)> [|t-a|(|\log|t-a||+C)]^{-1}, 
C=2\log n-\log
c''(H)+5, \ \text{for any critical value}\
a.\label{pb}\end{equation}
\end{corollary}
The Corollary follows from the previous Inequality and monotonicity
of the Poincar\'e metric.

\begin{proof} {\bf of (\ref{maxir}).} Let $\a_i$ be a path
from $t_0$ to $a=a_i$ from a marked path set, $|\a_i|\leq9$,
$\a=\la^{-1}\a_i$. As in the previous
Subsection (without loss of generality we consider that
the form $\omega$ has the type (\ref{omkl}) with
$\max|a_{kl}|=1$), one has
\begin{equation}|I(\hat t)|<n^4(2^{-2n}M_2)\leq M_2,\label{newnew}
\end{equation}
where  $M_2$ is the same as in (\ref{uppi3}) (recall that 
$\beta=\nu$). 
Let us estimate $M_2$: we show that 
\begin{equation}\log M_2<(1-\log c')R^6e^R.\label{m2<<}
\end{equation}
By elementary inequalities, the latter right-hand side is less than 
$\log M(R)$. This together with (\ref{newnew}) implies 
(\ref{maxir}). 

The linear combination $V$ of variations and $|t'|=|t|$ 
are estimated by inequalities (\ref{mr}) and (\ref{var<<}) 
respectively (proved below). 
Let us estimate the length of $\wt\a$: we show that 
\begin{equation}|\wt\a|<12R\log R.\label{mwta>}\end{equation}
By definition, the curve $\wt\a$ consists of the arcs of paths $\la$ 
and $\a_i$ lying in 
$\overline D_3\setminus\cup_iD_{\nu}(a_i)$. 
Those contained in $\a_i$ have total length less than 9, since 
$|\a_i|\leq9$. Those contained in $\la$ have total 
length no greater than 
$$|\la|_PM_3, \ M_3=(\min_{\wt\a}P(B))^{-1},\ |\la|_P\ \text{is 
the Poincar\'e length, thus,}\ |\la|_P\leq R+1.$$
Let us estimate $M_3$. 
Recall that the curve $\wt\a$, where the minimum in $M_3$ 
is taken, lies in 
$\overline D_3$ and its gap from the critical values is no less  
than $\nu$. 
This together with (\ref{pb}) and the inequality 
$|a|\leq2$ implies 
$$M_3\leq\max_{|t|\leq3}|t-a|(-\log\nu+C)\leq5(-\log\nu+C).$$
Inequality (\ref{3.5}) together with elementary inequalities implies 
that 
\begin{equation}C<\log R, \ -\log\nu<\log R.
\label{c>logr}\end{equation}
Therefore, $M_3<10\log R$. This together with the previous 
discussion and (\ref{3.5}) implies  that 
$$|\wt\a|<9+|\la|_P10\log R<11(R+1)\log R<12R\log R.$$
This proves (\ref{mwta>}). Substituting itself, (\ref{var<<}) and 
the inequality $|t'|=|t|<M_R$ (which follows from (\ref{mr})) 
 to the expression (\ref{uppi3}) of $M_2$ we get 
$$\log M_2<20n^{12}\frac{12R\log R+37e^RR\log R+5}{\nu}
-28n^4\log c'+2\log\frac{M_R}5.$$
By elementary inequalities and (\ref{3.5}), 
the latter right-hand side is less than 
$$e^RR^5\log R-R\log c'+12e^R\log R<(1-\log c')R^6e^R.$$
This proves (\ref{m2<<}) and (\ref{maxir}). 
\end{proof}

\begin{proof} {\bf of Theorem B1.} One has
\begin{equation}  \rho_R<5R. \label{rhr}\end{equation}
This follows from the fact that the diameter of $K_R=DP_R$ in the
Poincar\'e metric of $W$ is equal to $2R$ (by definition), and
the inequality
$$\frac{PU_R}{PW}|_{K_R}\leq\frac{e+1}{e-1}<\frac52.$$
The latter inequality is a particular case of the following more
general statement.

\begin{proposition} Let $W$ be a hyperbolic Riemann surface,
$U\subset W$ be a domain, $K\Subset U$ be a compact
set. Let $dist_{PW}(K,\partial U)\geq\sigma>0$. Then
$$\frac{PU}{PW}|_K\leq\frac{e^{\sigma}+1}{e^{\sigma}-1}.$$
\end{proposition}
\begin{proof} By monotonicity of the Poincar\'e metric as a function
of domain, it suffices to prove the Proposition in the case, when
$W=D_1$, $K=\{ 0\} ,$ $U$ is the Poincar\'e disc of radius
$\sigma$ centered at 0: in this case we prove the equality.
Indeed, let $r$ be the Euclidean radius of the latter disc. By
definition and conformal invariance of the Poincar\'e metric,
$$\frac{PU}{PD_1}(0)=r^{-1}. \ \text{One has}\ r^{-1}=
\frac{e^{\sigma}+1}{e^{\sigma}-1},$$
$$\text{since by definition,}\
\sigma=\int_0^r2\frac{ds}{1-s^2}=\log\frac{1+r}{1-r}.$$
This proves the Proposition.
\end{proof}

Let us estimate $B_{K_R,U_R}$. We show that
\begin{equation} B_{K_R,U_R}<(1-\log c'(H))e^{1.3R}.\label{bkr}
\end{equation}
Together with (\ref{rhr}) and (\ref{krur}), this implies Theorem
B1.

The set $K$ from the Main Lemma is contained in $K_R$, 
thus, $\log \max_{K_R}|I|\geq\log m$. Hence, by (\ref{maxir}),
\begin{equation} B_{K_R,U_R}<\log M(R)-\log m\label{burkr}
\end{equation}
We have shown at the end of 2.7 that
$$\log m\geq\log \Delta_0-(n^2-1)\log M_0-n^2\log n-\log2,$$
$$\Delta_0=(c'(H))^{6n^3}(c''(H))^{n^2}n^{-62n^3}, \
M_0=e^{\frac{2600n^{16}}{c''(H)}}(c'(H))^{-28n^4}.$$ 
Together with
(\ref{maxir}), (\ref{burkr}), (\ref{3.5}) and elementary 
inequalities this
implies (\ref{bkr}). Theorem B1 is proved modulo inequalities
(\ref{mr})-(\ref{var<<}).\end{proof}

\begin{proof} {\bf of (\ref{var<<}) modulo (\ref{mr}) and 
(\ref{betar}).} The expression 
$V=V_{\a,\nu}$ is a linear combination of variations of arguments 
along the pieces of the path $\a$ that lie either inside 
$\beta=\nu$- neighborhoods of the critical values of $H$, or 
outside $\overline D_3$. To estimate it from above, 
we use the following a priori upper bounds of variations.

Let $a$ be a critical value. By definition, 
 for any curve $l\subset B$ 
$$
\Var_l(t-a)=\int_l\frac{|dt|}{|t-a|}= 
\int_l{|dt|_P}\frac{(P(B))^{-1}}{|t-a|}\leq|l|_P\max_l 
\frac{(P(B))^{-1}}{|t-a|}$$ 
(here by $|l|_P$ we denote the Poincar\'e length). The 
latter ratio is estimated by (\ref{pb}): 
\begin{equation}\frac{(P(B))^{-1}}{|t-a|}<|\log|t-a||+C<7e^R\log R, 
\ \text{whenever}\ t\in\overline{DP_{R+1}}   
\label{pbcr}\end{equation} 
(the last inequality follows from (\ref{c>logr}) and (\ref{betar})). 
 Then by (\ref{pbcr}), 
\begin{equation}\Var_l(t-a)<7|l|_Pe^R\log R,\ \text{whenever} 
\ l\subset\overline{DP_{R+1}}.
\label{varn<}\end{equation}
Analogously,   for any critical value $a$ 
$$\Var_lt\leq|l|_P\max_l\frac{(P(B))^{-1}}{|t|}\leq 
|l|_P\max_l\frac{|t-a|}{|t|}\max_l\frac{(P(B))^{-1}}{|t-a|}.$$ 
Now let $l\subset\overline{DP_{R+1}}\setminus{\overline D_3}$. 
Then the former maximum in the previous  
right-hand side is no greater than $\frac53$, since $|t|>3$ on $l$ 
and $|a|\leq2$. Substituting this inequality and (\ref{pbcr}) 
to the same right-hand side yields 
\begin{equation}\Var_lt< 
\frac537|l|_Pe^R\log R<12|l|_Pe^R\log R, \ \text{whenever} \ 
l\subset\overline{DP_{R+1}}\setminus\overline{D_3}.
\label{beskvar}\end{equation} 

Let us estimate the expression $V=V_{\a,\nu}$. By definition, 
the variations in this expression are taken along 
the arcs of the  path $\a=\la^{-1}\a_i$ that lie 
either inside $D_{\nu}(a_j)$, or outside 
$\overline D_3$ (except for its final arc 
$\a(\tau',1]\subset D_{\nu}(a_i)$, $\a(\tau')\in
\partial D_{\nu}(a_i)$). 
By definition, the latter arc coincides with an arc 
of the path $\a_i$, and its complement in $\a_i$ is a curve 
lying in $\overline D_3$ 
outside the $\nu$- neighborhoods of the critical values 
(see 2.9). Therefore, the 
previous arcs, where the variations are taken, are disjoint from the 
path $\a_i$ and thus, are those of the path $\la$. The first sum 
in the expression of $V_{\a,\nu}$, 
which is $\nu$ times the sum of the variations along pieces 
of $\a$ near the 
critical values, is less than $7\nu(R+1)e^R\log R$. This follows 
from inequality (\ref{varn<}) applied to each piece and the 
inequality $|\la|_P\leq R+1$. Analogously, 
by the latter inequality and (\ref{beskvar}), the second sum in 
the expression of $V_{\a,\nu}$ is less than $36(R+1)e^R\log R$. 
The two previous upper bounds of the sums in $V_{\a,\nu}$ 
together with (\ref{3.5}) and inequality $\nu\leq\frac1{16}$ 
imply that  
$$V_{\a,\nu}<(36+7\nu)(R+1)\log Re^R<37R\log Re^R.$$
This proves (\ref{var<<}). 
\end{proof}

\begin{proof} {\bf of (\ref{mr}).} Let $a$ be a critical value of 
$H$, $t\in\overline{DP_{R+1}}$. Let us prove that $|t|<M_R$: 
this will imply (\ref{mr}).
It follows from definition and (\ref{pb}), (\ref{c>logr}) that 
\begin{equation}R+1\geq\int_{|t_0-a|}^{|t-a|}\frac{|ds|}{s(|\log s|
+C)}, \text { where } C<\log R.
\label{rgeq}\end{equation}

 By definition, $|a|\leq2$,
$|t_0|\leq3$, so, $|t_0-a|\leq5$. Suppose $|t|>7$ (if not, then 
the inequality $|t|<M_R$ follows immediately, since 
$M_R>7$ (by (\ref{3.5}) and elementary inequalities)). 
Hence, $|t-a|>5$. Put
$u=\log s$. Then the latter integral is greater than
$$\int_5^{|t-a|}\frac{ds}{s(\log s+C)}=
\log(u+C)|_{\log 5}^{\log|t-a|}$$ By elementary inequalities, the
latter right-hand side is greater than
$$\log\log|t-a|-\log(C + 2).$$
This together with (\ref{rgeq}) implies that 
$$\log|t-a|<e^{R+1}(C+2).$$
This together with inequality $|a|\leq2$, (\ref{c>logr}) 
and elementary inequalities implies (\ref{mr}).
\end{proof}

\begin{proof} {\bf of (\ref{betar}).} It suffices to show that
for any critical value $a$
\begin{equation} |t-a|>\beta_R\ \text{for any} \ t\in 
\overline{DP_{R+1}}.
\label{tabr}
\end{equation}
It follows from formula for $\beta_R$ in (\ref{betar}), 
inequality (\ref{3.5}), choice of $t_0$ and elementary inequalities 
that    
\begin{equation} \beta_R<\nu=\frac{c''(H)}{4n^2}\leq|t_0-a|.
\label{bnu}\end{equation}
Thus, if $|t-a|\geq\nu$, then inequality (\ref{tabr}) holds.
Let us prove (\ref{tabr}) assuming that $|t-a|<\nu$. To do this,
we use the fact that under this assumption the integral in
(\ref{rgeq}) is greater than
$$\int_{\nu}^{|t-a|}\frac{ds}{s(|\log s|+ C)}=
\log(u+ C)|_{-\log\nu}^{-\log|t-a|}$$
$$>\log\log(|t-a|^{-1})-\log(-\log\nu+ C).$$
This together with (\ref{rgeq}), (\ref{c>logr}) 
and elementary inequalities implies (\ref{tabr}).
Inequality (\ref{betar}) is proved.
\end{proof}

\subsection{Proof of the Geometric Lemma \ref{l1.3}}


We prove Lemma \ref{l1.3} by induction  in $s$. For $s=1$ it is a
direct consequence of the definition of vanishing cycle. Indeed,
in this case $[\gamma]=[\gto]$ is the cycle vanishing along the
segment $[t_0, a_1]$.

Let the statement of Lemma \ref{l1.3} be proved for all $s<N$. Let
us prove it for $s=N$.

Denote  by $A_i$, $i=1,\dots,N$ the critical points located inside
$\gamma$. Let $a_i$ be the corresponding critical values. Without
loss of generality suppose that the value $H(x,y)$ decreases
locally when the point $(x,y)$ in the real plane moves from the
oval $\gto$ inside the domain bounded by this oval (this may be
achieved by changing the sign of $H$). There is a critical point
$A_i$ such that the oval $\gto$ extends up to a continuous family
of real ovals
\def\gt{\gamma_t}
$\gt\subset S_t$ on the semiinterval $(a_i,t_0]$ so that the limit
$\lim_{t\to a_i}\gt$ is a loop with the base point $A_N$ (see Fig.
5a, b). This loop is a connected component of a critical level
that contains only one singular point of $H,$ because  $H$ is
ultra-Morse. Hence, the limit loop may be either an eight-shaped
figure, or a simple loop, see Figures 5a and 5b respectively.
Geometric Lemma is proved below in case of the eight-shaped
figure, which is a union of two simple loops  $\Gamma_1$ and
$\Gamma_2$ that are  disjoint (outside $A_N$) and bound disjoint
domains (see Fig. 5a). Another case depicted at Fig. 5b is treated
analogously.

Choosing appropriate numeration of the $a_j$'s, suppose that
$i=N$. Without loss of generality we may assume that $a_N=0$,
$A_N=0$ (this may be achieved by  real translations in the source
and target of the map $H$).

When $t\in\mathbb R_+$ passes  through 0 to $\mathbb R_-$, the
loop $\Gamma_1 \Gamma_2 $ generates a pair of ovals $\gt^i$,
$i=1,2$, in the real level curve $H(x,y)=t$: the oval $\gt^i$ lies
in the domain bounded by the curve $\Gamma_i$ and tends to
$\Gamma_i$, as $t\to0$ (see Fig.5a). Suppose that the  curves
$\gto$, $\Gamma_i$ and $\gamma_{t}^i$ are oriented
counterclockwise. All the critical points $A_j, j < N,$ are
contained in the domains bounded by the ovals $\gt^i$. By the
induction assumption, each oval $\gt^i$ satisfies the statement of
Lemma \ref{l1.3}: (\ref{1.5}) holds for $\gto$ replaced by
$\gamma_t^i$, $s = N$ replaced by $s<N$.

\begin{figure}[ht]
  \begin{center}
   \epsfig{file=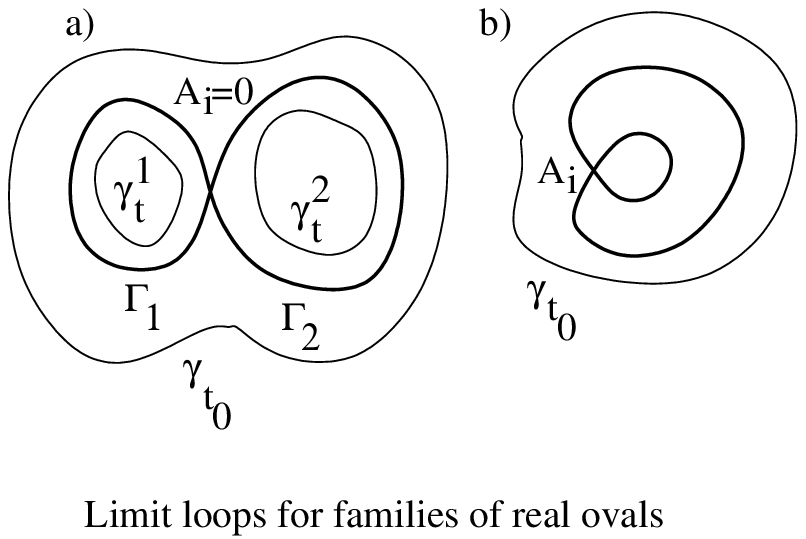}
    \caption{}
    \label{fig:5}
  \end{center}
\end{figure}

We proceed below the induction step for the case when $t_0$ is
small.
In the general case we connect $t_0$ to a small $t_0'\in\mathbb
R_+$ by a segment $\a'=[t_0,t_0']$. The corresponding family of
ovals  $\gamma_t$ starting with $\gto$ is well-defined over $\a'$.
If $\a'$ does not contain critical values of $H$, then the
decomposition (\ref{1.5}) of $\gamma_{t_0'}$ with $\delta_j$
replaced by the cycles vanishing along $\a'\circ\a_j$ (this
decomposition for $\gamma_{t_0'}$ is proved below) extends along
$\a'$ to decomposition (\ref{1.5}) of $\gto.$ Otherwise, we
replace $\a'$ by its small deformation to the upper half-plane.
(Since $\gamma_t$ is well-defined on the segment, the result of
its extension to $t_0$ along the deformed path $\a'$ is not
changed: it is the real oval $\gamma_{t_0}$.)

Let us now prove (\ref{1.5}) for $t_0$ small. Consider the
semicircular path $\tau(\theta)=t_0e^{i\theta}$,
$\theta\in[0,\pi]$, which goes around the zero critical value in
the upper half-plane. Let $\delta_N\in H_1(S_{t_0},\mathbb Z)$,
$\delta_N'\in H_1(S_{-t_0},\mathbb Z)$ be the cycles vanishing
along the real segments going from $\pm t_0$ to 0. By definition,
the cycle $\delta_N'$ is obtained as the extension of the cycle
$\delta_N$ along the path $\tau(\theta)$ . Let us show that
the curve $\gto$ admits a homotopy by curves $\gt$
in complex level lines $H(x,y)=t$ along the path $t=\tau(\theta)$
so that
\begin{equation}
[\gamma_{-t_0}]=[\gamma_{-t_0}^1]+ [\gamma_{-t_0}^2]\pm [
\delta_N']. \label{1.25}
\end{equation}
Together with the decompositions (\ref{1.5}) for $\gamma_{-t_0}^i$
that is valid by the induction hypothesis, this implies
(\ref{1.5}) for $\gto$ and  completes the induction step.

To construct the homotopy mentioned above of the oval $\gto$, let us
consider the real local analytic coordinate system $(x',y')$ in a
neighborhood of zero critical point such that $H=x'y'$; it exists
by the Morse lemma. {\it Then locally near 0, the curves
$\Gamma_i$} are intervals in the new coordinate lines. By the
choice of $t_0$, we may suppose that the previous neighborhood
contains the square $U$ centered at 0 (in the new coordinates
$(x',y')$) whose sides are parallel to the coordinate axes and
have length $2\sqrt{t_0}$. The curve $U \cap \gto$ lies  in the
first and third quadrants of this chart: $x',y'>0$; $x',y'<0$. The
entire curve $\gto$  is split by the points
$b_+=(\sqrt{t_0},\sqrt{t_0})$, $b_-=(-\sqrt{t_0},-\sqrt{t_0})$
into two arcs denoted by $\Gamma_i(t_0)$, $i=1,2$. Suppose  that
the intersection of the domain bounded by $\Gamma_1$ with $U$
belongs to  the quadrant $x>0,y<0$. Then
$(0,-\sqrt{t_0}),(\sqrt{t_0},0)\in\Gamma_1$, and the curve
$\Gamma_1(t_0)$ is oriented from $b_-$ to $b_+$; $(-\sqrt{t_0},0),
(0,\sqrt{t_0})\in\Gamma_2$, and the curve $\Gamma_2(t_0)$ is
oriented from $b_+$ to $b_-$ (see Fig.6a). The vanishing cycle
$\delta$ is represented by the circle
$\tilde{\delta}=\{(-\sqrt{t_0}e^{i\psi},\sqrt{t_0}e^{-i\psi}) \ |
\ \psi\in[0,2\pi]\}$.

Our goal is to construct the family $\gamma_{\tau(\theta)}$ and
then check  (\ref{1.25}).  Below we construct the homotopy 
$\Gamma_i(\tau(\theta))$ of each arc $\Gamma_i(t_0)$ along the
path $\tau(\theta)$ (as a family of arcs in complex level curves
$H=\tau(\theta)$) so that

1) $\Gamma_i(-t_0)=\gamma_{-t_0}^i$;

2) the arc $\Gamma_1(\tau(\theta))$ starts at
$(-\sqrt{t_0}e^{i\theta},-\sqrt{t_0})$ and
ends at $(\sqrt{t_0},\sqrt{t_0}e^{i\theta})$;

3) the arc $\Gamma_2(\tau(\theta))$ starts at
$(\sqrt{t_0}e^{i\theta},\sqrt{t_0})$ and
ends at $(-\sqrt{t_0},-\sqrt{t_0}e^{i\theta})$.

Then we put
$$\Gamma_{\pm}(\theta)=\{(\pm\sqrt{t_0}e^{i\phi},
\pm\sqrt{t_0}e^{i(\theta-\phi)})|\  0\leq\phi\leq\theta\} \
\text{(with a natural orientation)},$$

and a representative of the class $[\gamma_{\tau(\theta)}]$ may be
constructed as a product of four curves:
$$\gamma_{\tau(\theta)}=\Gamma_1(\tau(\theta))\Gamma_+(\theta)\Gamma_2(\tau(\theta))
\Gamma_-(\theta).$$

\begin{figure}[ht]
  \begin{center}
   \epsfig{file=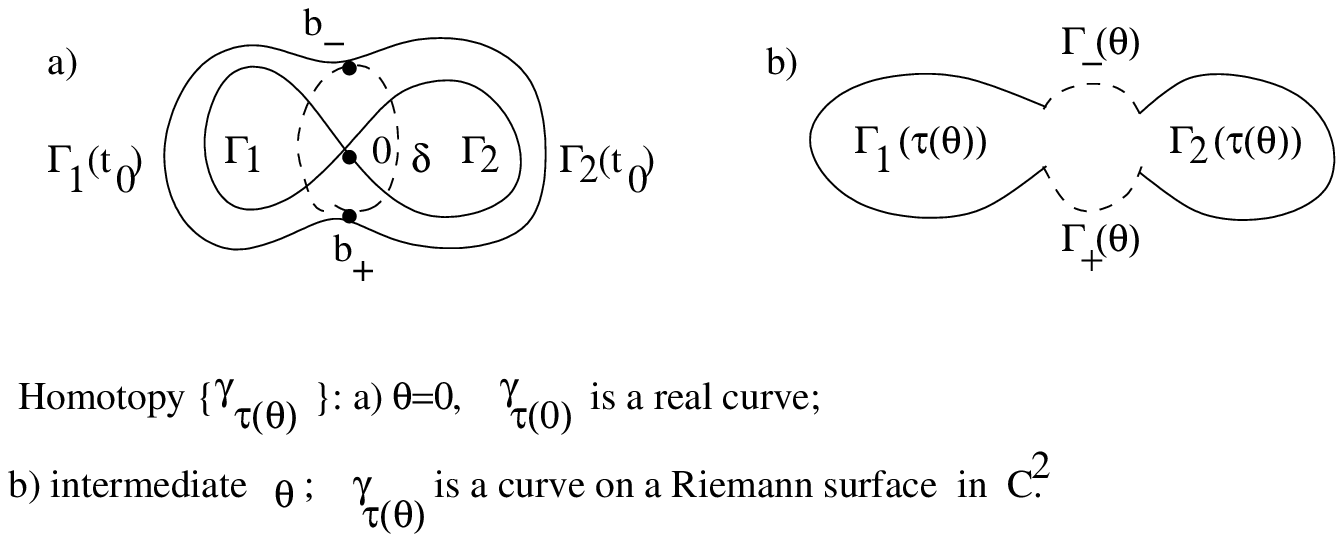}
    \caption{}
    \label{fig:6}
  \end{center}
\end{figure}

By construction, the curves $\gamma_{\tau(\theta)}\subset\{
H(x,y)=\tau(\theta)\}$ are closed (with a well-defined
orientation) and depend continuously on $\theta$,
$\gamma_{\tau(0)}=\gto.$ Note that $ [ \delta_N']=
[\Gamma_+(\pi)\Gamma_-(\pi)].$  This implies (\ref{1.25}) and
proves Lemma \ref{l1.3} modulo existence of families
$\Gamma_i(\tau(\theta))$ satisfying assumptions 1)-3) above.

We construct the family $\Gamma_1(\tau(\theta))$ only, the family
$\Gamma_2(\tau(\theta))$ is constructed analogously. To  do this,
consider an arbitrary increasing parametrization
$\psi:u\in[0,1]\mapsto\Gamma_1$ of the oriented curve $\Gamma_1$,
$\psi'\neq0$, $\psi(0)=\psi(1)=0$. Recall that the local starting
branch (going from 0) of $\Gamma_1$ lies in the negative $y'$-
semiaxis, so, the coordinate $-y'$ increases locally along this
branch. Analogously, the coordinate $-x'$ increases locally along
the final branch entering 0 of $\Gamma_1$. Let us choose the
previous parametrization so that
$$u=-y'\circ\psi(u)\ \text{near}\ u=0; \ u=-x'\circ\psi(u)+1
\ \text{near}\ u=1.$$
The mapping $\psi$ extends up to a locally invertible
$C^{\infty}$ mapping
of $\mathbb C^2$- complex neighborhood of the initial parameter segment $[0,1]$
(we consider that the coordinates in the new parameter space $\mathbb C^2$  are $(u,v)$,
and the previous segment $[0,1]$ lies in the complex $u$- axis).
One can choose the previous extension of $\psi$ so that

a) the previous equalities hold in complex
neighborhoods of the points $(0,0)$, $(1,0)$:

\begin{equation}
u=-y'\circ\psi(u,v)\ \text{near} \ (0,0);\ u=-x'\circ\psi(u,v)+1
 \ \text{near}\ (1,0);\label{1.26}\end{equation}

b) the level curves of the pull-back $H\circ\psi$ of $H$
(except for the lines
$\{u=0,1\}\subset\{ H\circ\psi=0\}$) are transversal
to the lines $u=const$.

Then $(u,H)$ are well defined global coordinates on the complement of the parameter
domain to the latter pair of lines $u=0,1$.
Let $\widetilde{\Gamma_1(t_0)}$ be the lifting
to the parameter domain of the arc $\Gamma_1(t_0)$,
$\widetilde{\Gamma_1(\tau(\theta))}$ be its image under the mapping preserving
the coordinate $u$ and
multiplying the coordinate $H\circ\psi$ by $e^{i\theta}$. The arc
$\Gamma_1(\tau(\theta))=\psi(\widetilde{\Gamma_1(\tau(\theta))})$ is the
one we are looking for. Indeed, it lies in the complex level curve
$H=\tau(\theta)$ by construction. It starts at $(-\sqrt{t_0}e^{i\theta},
-\sqrt{t_0})$ and ends at $(\sqrt{t_0},\sqrt{t_0}e^{i\theta})$ by (\ref{1.26}) and
the equality $H=x'y'$. It follows from construction and (\ref{1.26}) that
$\Gamma_1(-t_0)=\gamma_{-t_0}^1$. Lemma \ref{l1.3} is proved.

\section{Estimates of the number of zeros of Abelian integrals near
the critical values}

In this section we prove Theorem A2, see 1.8, and Theorem B2,
stated below. Together with Theorem A1 (whose proof is completed
in section 2) Theorem A2 implies Theorem A. Together with Theorem
B1 ( whose proof is completed in section 3), Theorem B2 implies
Theorem B.

We have three statements to discuss:

1. Theorem A2 in the case when the endpoints of the interval
 considered are all finite;

2. Theorem A2 in the case when one of these endpoints is infinite;

3. Theorem B2.

These statements will be referred to as cases 1,2,3 below.

It appears that cases 1 and 3 are very close to each other.

\subsection{Argument principle, KRY theorem and
Petrov's method}

All the three  cases are treated in a similar way. We want to
apply the argument principle.

The estimates near infinity are based on the argument principle
only. The estimates near finite critical points use the Petrov's
method that may be considered as a generalization of the argument
principle for multivalued functions. The increment of the argument
is estimated through the Bernstein index of the integral, bounded
from above in the previous sections. The relation between these
two quantities is the subject of the Khovanskii-Roitman-Yakovenko
(KRY) theorem and Theorems 4.3, 4.4 stated below. It seems
surprising that these theorems were not discovered in the
classical period of the development of complex analysis. The
preliminary version of the theorem was proved in \cite{31}, the
final one in \cite{16}. One of the two inequalities in this
theorem is proved by the second author (Yu.S.Ilyashenko,
\cite{KRY}) in a stronger form (Theorem \ref{t6.2} below): an
explicit formula for the constant in the estimate is written, and
$\Bbb C$ is replaced by an arbitrary Riemann surface.

At this spot we begin the proof of Theorem A2 in case 1. Recall
the statement of the theorem in case 1.

\medskip

{\bf Theorem A2 (Case 1).} {\it Let $ a \not = \infty, b \not =
\infty,.$ Then}
$$
\#\{t \in (a,l(t_0))\cup(r(t_0),b)  \mid I(t) = 0\}< (1-\log
c')e^{\frac{4800}{c''}n^4},
$$
where $l(t_0)$ and $r(t_0)$ are the same as at the beginning of
1.7.
\medskip

We will prove that
\begin{equation}
\#\{t \in (a,l(t_0))  \mid I(t) = 0\}< \frac 12 (1-\log
c')e^{\frac{4800}{c''}n^4}. \label{6.1}
\end{equation}
Similar estimate for $(r(t_0),b) $ is proved in the same way.
These two estimates imply Theorem A2.

Let $\Pi=\Pi(a)$ be the same as in (\ref{1.2b}).
\begin{lemma} \label{l6.1}
Inequality (\ref{6.1}) holds provided that in
(\ref{6.1}) the interval $(a,l(t_0))$ is replaced by $\Pi.$
\end{lemma}

Lemma \ref{l6.1} implies (\ref{6.1}) because $(a,l(t_0)) \subset \Pi.$ Let
$${\Pi}_{\psi}=\{ t \in \Pi \mid \psi \le |t-a| \le \nu \}$$

\begin{lemma} \label{l6.2} Lemma \ref{l6.1}
holds provided that in (\ref{6.1}) the
domain $\Pi$ is replaced by $\Pi_{\psi}.$
\end{lemma}

Lemma \ref{l6.2} implies Lemma \ref{l6.1}, because

$$\Pi = \cup_{\psi > 0}\Pi_{\psi}.$$

\begin{proof} {\bf of Lemma \ref{l6.2}.} The proof of
this lemma occupies this and
the next four subsections. We have

$$
\partial {\Pi}_{\psi}= \Gamma_1 \Gamma_2 \Gamma_3 \Gamma_4.
$$

As sets, the curves $ \Gamma_j$ are defined by the formulas below;
the orientation is defined separately:

$$\Gamma_1=\{t \mid |t-a|=\nu, |arg(t-a)| \le 2 \pi \}=\Gamma_a$$
$$\Gamma_3=\{ t \mid |t-a|= \psi , |arg(t-a)| \le 2\pi \} $$
$$\Gamma_{2,4}=\{ t \mid \psi \le |t-a| \le \nu, arg(t-a)=\pm2\pi\}.$$

The curve $ \Gamma_1 $ is oriented counterclockwise, $\Gamma_2$ is
oriented from the right to the left, $ \Gamma_3 $ is oriented
clockwise, $ \Gamma_4 $ is oriented from the left to the right.

Let $\#\{t \in (a+\psi,l(t_0))  \mid I(t) = 0\} = Z_{\psi}.$ Denote by
$R_{\Gamma}(f)$ the increment of the argument of a holomorphic function
$f$ along a curve $\Gamma$ ($R$ of Rouchet). Recall that $V_{\Gamma}(f)$ 
denotes the variation of the argument of $f$ along $\Gamma$.  
Obviously, $\mid R_{\Gamma}(f)\mid  \le V_{\Gamma}(f).$
\end{proof}

In assumption that $I \not = 0 $ on $ \partial {\Pi}_{\psi},$ the
argument principle implies that

\begin{equation}Z_{\psi} \le \frac 1{2\pi} R_{\partial {\Pi}_{\psi}} (I) \le \frac 1{2\pi}
\sum_1^4 R_{\Gamma_j} (I). \label{6.2}\end{equation}

The first term in this sum is estimated by the modified KRY
theorem, the second and the forth one by the Petrov method, the
third one by the Mardesic theorem. The case when the above
assumption fails is treated in 4.3.

\subsection{Bernstein index and variation of argument}

The first step in  establishing a relation between variation of
argument and the Benstein index was done by the following KRY
theorem.

Let $U$ be a connected and simply connected domain in $\Bbb C,$
$\Gamma\subset U$ be a (nonoriented) curve, $f$ be a bounded
holomorphic function on $U$.
\medskip

{\bf KRY theorem, \cite{16}} {\it For any tuple $U,\Gamma \subset U $
as above and a  compact set $K\subset U$ there exists a geometric
constant $\alpha=\alpha (U,K,\Gamma),$ such that}
$$ V_{\Gamma}(f)\le \alpha B_{K,U}(f).$$
\medskip

In \cite{16} an upper estimate of the    Bernstein index through the
variation of the argument along $ \Gamma = \partial U$ is given;
we do not use this estimate. On the contrary, we need an improved
version of the previous theorem with $\alpha$ explicitly written
and $ U$ being a domain on  a Riemann surface. These two goals are
achieved in the following two theorems.

Let $|\Gamma|$ be the length, and $ \kappa (\Gamma ) $ be the
total curvature of a curve on a surface endowed with a Riemann
metric.

\begin{theorem} \label{t6.1} Let
$\Gamma \subset U'' \subset U' \subset
U \subset \Bbb C$ be respectively a curve, and three open sets in
$\Bbb C.$ Let $f: U \to \Bbb C$ be a bounded holomorphic function,
$f|_{\Gamma} \not = 0.$ Let $\e < \frac {1}{2}$ and the following
gap conditions hold:
\begin{equation}
\rho (\Gamma , \partial U'') \ge \e , \ \rho (U'', \partial U')
\ge \e , \ \rho (U',\partial U) \ge \e . \label{1}\end{equation}

Let $D > 1$ and the following diameter conditions hold:
\begin{equation}
\text { diam }_{\text {int}}U'' \le D, \ \text { diam }_{\text
{int}}U' \le D.  \label{2}
\end{equation}
Then
\begin{equation}
V_\Gamma (f) \le B_{U'',U}(f)(\frac{\mid \Gamma \mid}{\e}
+ \kappa (\Gamma )
+ 1)e^{\frac {5D}{\e }}.   \label{3}
\end{equation}
\end{theorem}
\begin{theorem} \label{t6.2} Let $\Gamma \subset U'' \subset U' \subset
U \subset W$ be respectively a curve, and three open sets in a
Riemann surface $W.$  Let $f: U \to \Bbb C$ be a bounded
holomorphic function, $f|_{\Gamma} \not = 0.$ Let $\pi: W \to \Bbb
C$ be a projection which is locally biholomorphic, and the metric
on $W$ is a pullback of the Euclidean metric in $\Bbb C.$
Let $\e < \frac {1}{2}$ and the following gap conditions hold:
\begin{equation}
 \pi \text {-gap }(\Gamma , U'') \ge \e , \ \pi \text {-gap
}(U'',U') \ge \e , \ \pi \text {-gap }(U',U) \ge \e .   \label{4}
\end{equation}
Let $D > 1$ and the following diameter conditions hold:
\begin{equation}
D > 1, \  {\text {diam }}_{int}U''\le D, \ {\text {diam
}}_{int}U'\le D    \label{5}
\end{equation}
Then inequality (\ref{3}) holds.
\end{theorem}
These theorems are proved in \cite{KRY}.


Recall that intrinsic diameter and $\pi $-gap are defined in 1.6.

We can now estimate from above the first term in the sum (\ref{6.2}).
The estimate works in both cases when $a$ is finite or infinite.
Let for simplicity, as in (\ref{1.3a}),
$$
A = e^{\frac {n^4}{c''}}.
$$

\begin{lemma} \label{l6.3}  Let $H$ be a
balanced polynomial of degree $n+1\geq3$. Let I be the same
integral as in (1.1). Let K be a compact set mentioned in the Main
Lemma, and $\Gamma_1 = {\Gamma}_a$ be the same as in this Lemma
($a$ may be infinite). Then
\begin{equation}
V_{{\Gamma}_1}(I)< (1 - \log c') A^{4700}.    \label{6.3}
\end{equation}
\end{lemma}


\begin{proof} The lemma follows immediately from Theorem \ref{t6.2} and
the Main Lemma. To apply Theorem \ref{t6.2}, let us take $I$ for $f,$ the
universal cover over $B$ for
$W$ with the natural projection $\pi :
W \to \Bbb C$ and metric induced from $ \Bbb C$ by this
projection. This metric on $W$ is called Euclidean. Let $K$ and
$U$ be the same as in the Main Lemma. Take this $U$ for the domain
$U$ to apply Theorem \ref{t6.2}. Recall that $U$ is the minimal simply
connected domain that contains the $\frac {\nu }{2}$-\nbd of $K$ in
$U$ in the Euclidean metric on $W, \ \nu $ is the same as in
(\ref{1.2}). Let $\e = \frac {\nu }{6},$ that is
\begin{equation}
\e = \frac {c''}{24n^2}.  \label{6.3'}
\end{equation}
Let $U''$ and $U'$ be the minimal simply connected domains in $W$
that contain $\e $-\nbd of $K$ and $U''$ respectively. Note that
$\Gamma_1 = \Gamma \subset K.$ Then gap condition (\ref{4}) with $\e $
from (\ref{6.3'}) holds. Moreover, $diam_{int} U\le diam_{int}
U^{\prime}+2\varepsilon \le diam_{int} K+4\varepsilon .$ Hence,
diameter condition (\ref{5}) holds with $D<38n^2$ by (\ref{1.6}).
Thus

$$
e^{\frac {5D}{\e}} < A^c, \text { where } A = e^{\frac
{n^4}{c''}}, \ c = 5\times38\times24 < 4600.
$$

This factor $A^c$ is the largest one in the estimate for
$V_{\Gamma_1} (I).$

By inequality (\ref{1.8}) from the Main Lemma, $B_{K,U} \le (1 - \log c')A^2.$ 
By the monotonicity of Bernstein index (that follows directly
from the definition), $B_{U'',U} < B_{K,U}.$ At last,
$$
\frac{|\Gamma_1|}{\e} + \kappa (\Gamma_1) + 1=
24\pi+ 4 \pi + 1 << A.
$$

Now, inequality (\ref{3}) proves the  lemma.
\end{proof}

The Corollary below is used in the next subsection.

\begin{remark} Lemma \ref{l6.3} 
remains valid if in its hypothesis the integral $I$ is replaced by
an integral $J$ over the cycle vanishing at the critical value $a$
of $H$. The proof of this modified version of Lemma 6.3 
repeats that of the original one with the following change:
we use the Modified Main Lemma instead of the Main Lemma.
\end{remark}

\begin{corollary} \label{c6.1} Suppose that the integral $J$ with a
real integrand $\omega$ is taken over a local vanishing cycle
$\delta_t$ corresponding to the real critical value $a$. Then the
number of zeros of $J $ in the disk centered at $a$ of radius
$\nu=\frac{c''}{4n^2}$ admits  the following upper estimate:
\begin{equation}
N_J:= \#\{t \in \Bbb C \mid  |t-a|<\nu, J(t)=0\} \le \frac
1{2\pi}(1-\log c')A^{4700}  \label{6.3a}
\end{equation}
\end{corollary}

This follows from the modified Lemma \ref{l6.3} and the argument
principle.

\subsection{Application of the Petrov's method}

The Petrov's method applied below is based on a remark that the
magnitude of the increment of the argument of a nonzero function
along an oriented curve is no greater than the number of zeros of
the imaginary part of this function increased by 1 and multiplied  by
$\pi .$ Indeed, at any half circuit around zero, a planar curve
crosses an imaginary axis at least once. The method works when the
imaginary part of a function appears to be more simple than the
function itself.

Let $\delta_t \in H_1(t)$ be the local vanishing cycle at the
point $a.$ Let $\omega $ be the same real form as in integral
(1.1). Let $J$ be the germ of integral $J(t) = \int_{\delta_t}
\omega $ along the cycle $ \delta_t$, which is a local vanishing
cycle at $t = a.$ Note that $J$ is single-valued in any
simply connected \nbd of $a$
that contains no other critical values of $H.$ Let $l_0 =
(\gamma_t, \delta_t)$ be the intersection index of the cycles
$\gamma_t$ and $\delta_t.$ As the cycle $\gamma_t$ is real and $H$
is ultra-Morse, $l_0$ may take values $\pm 1, \pm 2$ only
(Lemma \ref{l2.4}).
Let
$$
\Gamma_0 =\{ t \in\mathbb R\mid te^{2\pi i}\in \Gamma_2\} .
$$
Then by the Picard-Lefschetz theorem
$$
I\mid_{\Gamma_2} = (I + l_0J)\mid_{\Gamma_0}, \ I\mid_{\Gamma_4} =
(I - l_0J)\mid_{\Gamma_0}.
$$
\begin{proposition} \label{p6.1} The integral J is purely imaginary  on
the real interval $(a,b)$.
\end{proposition}
\begin{proof} Recall that the form $\omega$ and the polynomial H
are real. Then
$$J(t)=-\overline{J(\overline{t})}.$$ Indeed, $\omega=Q(x,y)dx$. The
involution ${\mathbf i}:(x,y)\mapsto (\overline x,\overline y)$
brings the integral $J(t)=\int\limits_{\delta_t}Qdx$ to
$\int\limits_{{\mathbf i}\delta_t}\overline Q d \overline x=
\int\limits_{-\delta_{\overline t}}\overline Q d \overline x=
-\overline{\int_{\delta_{\bar t}}Qdx}=-\overline{J(\overline t)}.$
On the other hand, for real $t$ we have $t=\overline t$ and
$\delta_{\overline t}=\delta_t$. Hence, $J(t)=-\overline {J(t)}$
for $t \in (a,b).$ This implies Proposition 6.1.
\end{proof}
\begin{corollary} \label{c6.2} 
$$
Im I \mid_{ \Gamma_{2,4}}=\pm l_0J \mid_{ \sigma}.
$$
\end{corollary}

\begin{proof} This follows from Proposition \ref{p6.1}, 
Picard-Lefschetz theorem and the reality    of $I$ on   $\sigma.$
\end{proof}

Suppose now that $I$ has no zeros on $\Gamma_2$ and $\Gamma_4.$
Then
\begin{equation}
\left| R_{\Gamma_{2,4}}(I)\right| \le \pi ( 1 + N),  \text { where
} N = \# \{ t \in \Gamma_0\mid J(t) = 0\}. \label{6.4}
\end{equation}
Obviously, $N \le N_J,$ see (\ref{6.3a}). The right hand side of this
inequality is already estimated from above in Corollary \ref{c6.1}.
Hence,
$$
\left| R_{\Gamma_{2,4}}(I)\right| \le \pi (1 - \log c')A^{4700}.
$$

Suppose now that $I$ has zeros on $\Gamma_2$ (hence on $\Gamma_4$,
by Proposition \ref{p6.1}). Indeed, its real part is the same at the
corresponding points of $\Gamma_2, \Gamma_0, \Gamma_4, $ and the
imaginary  parts of $I|_{\Gamma_2}$ and $ I|_{\Gamma_4}$ are
opposite at the corresponding points.) In this case we replace the
domain $\Pi_\psi $ by $\Pi_\psi ' $ defined as follows.

The curves $\Gamma_{2,4}$  should be modified. A small segment of
$\Gamma_{2}$ centered at zero point of $I$ that contains no other
zeros of $J,$ should be replaced by an upper half-circle having
this segment as a diameter and containing no zeros of $J.$ A
similar modification should be done for $\Gamma_4 $ making use of
lower half-circles. Denote the modified curves by $
\Gamma_{2,4}'.$ Let $\Pi_\psi'$ be the domain bonded by the curve
\begin{equation}
\partial {\Pi}_{\psi}'= \Gamma_1 \Gamma_2' \Gamma_3 \Gamma_4'.
\label{6.4a}
\end{equation}
It contains $\Pi_{\psi}$, and we will estimate from above the
number of zeros of $I$ in $\Pi_\psi ' $ still using the argument
principle. The increment of $\arg I$ along $ \Gamma_1 $ is already
estimated in 4.2. Here we give an upper bound for the increment of
$\arg I$ along $\Gamma_{2,4}'.$ The increment along $\Gamma_3$ is
estimated in the next subsection.

\begin{proposition} \label{p6.2} Let $N$ be the same as in (\ref{6.4}).
Then
\begin{equation}
\mid R_{\Gamma_2'} (I)\mid  \le \pi (2 N + 1). \label{6.5}
\end{equation}
\end{proposition}

\begin{proof} Let $I$ have zeros $b_j \in \Gamma_2, \ j = 1,...,k,
$ the number of occurrence  of $b_j$ in this list equals its
multiplicity. Note that
\begin{equation}
\text {Im } I|_{\Gamma_2} = \pm l_0J\circ \pi,\
l_0=<\delta,\gamma>\neq0.   \label{6.5b}
\end{equation}
Hence, at the points $b_j, \  J $ has zeros of no less
multiplicity than $I.$  Hence, the total multiplicity $k'$ of zeros
of $J$ at the points $b_j \in \Gamma_2, \ j = 1,...,k, $ is no less
than $k.$ Let $J$ have $s$ zeros on $\Gamma_2'. $ We have: $k' \ge
k, \ s \le N-k' \le N-k.$ Let $\sigma_1, ... , \sigma_q, \ q \le k
+ 1,$ be the open intervals into which the curve $ \Gamma_2$ is
divided by the points $b_j.$ Let $s_j$ be the number of zeros of
$J$ on $\sigma_j, \ \sum_1^q s_j = s.$ Let
$$
R_j = R_{\sigma_j} (I).
$$
Then
$$
R_j \le \pi (s_j + 1).
$$
Hence,
\begin{equation}
\mid R_{\Gamma_2'}( I)\mid  \le \pi ( k + \sum_1^q ( s_j + 1 ))
\le \pi (2 k +1+ s) \le \pi (2 k'+1 +  s) \le \pi (2N+1).
 \label{6.5a}
\end{equation}
\end{proof}


\subsection{Application of the Mardesic theorem}

\begin{proposition} \label{p6.3} Let $I$ be the integral (1.1), and
$\Gamma_3$ be the same as in 4.1. Then for $\psi $ small enough,
\begin{equation}
\left| R_{\Gamma_3} (I)\right| \le \pi (4 n^4 +1). \label{6.6}
\end{equation}
\end{proposition}
\begin{proof} Let $J$ and $l_0$ be the same as in the previous
subsection. Let $a = 0,$ and $I(e^{2\pi i}t)$ means the result of the
analytic extension of $I$ from a value $I(t)$ along a curve
$e^{2\pi \ph }t, \ \ph \in [0,1].$ By the Picard-Lefshetz theorem,
for small $t$
$$
I(e^{2\pi i}t) = I(t)+ l_0J(t).
$$
Consider the function
$$
Y(t) = I(t)- l_0\frac {\log t}{2\pi i}J(t).
$$
This function is single-valued because the increments of both terms
$I$ and $Y$ under the analytic extension over a circle centered at
$0$ cancel. The function $I$ is bounded along any segment ending
at zero, and $J$ is holomorphic at zero, with $J(0) = 0.$ Hence, $Y$
is holomorphic and grows no faster than $\log |t|$ in a punctured
\nbd of zero.
(In fact, it is bounded in the latter neighborhood:
$|J(t)\log t|\leq c|t||\log t|\to0$, as $t\to0$.)
By the removable singularity theorem, it is
holomorphic at zero. Hence,
\begin{equation}
I(t) = Y(t) + l_0\frac {\log t}{2\pi i}J(t)   \label{6.7}
\end{equation}
with $Y$ and $J$ holomorphic. We claim that {\it the increment of
the argument of $I$ along $\Gamma_3$ for $\psi $ small is bounded
from above through $\text {ord}_0 J,$ the order of zero of $J$ at
zero.} The latter order is estimated from above by the following
theorem by Mardesic:

\begin{theorem} \cite{17}. The multiplicity of any zero of the
integral $I$ (or $J$) taken at a point where the integral is
holomorphic does not exceed $n^4.$
\end{theorem}

The function (\ref{6.7}) is multivalued. The proof of the latter
claim is based on the following simple remark. Let $ f_1, f_2$ be two
continuous functions on a segment
$ \sigma\subset\mathbb R$, and $|f_1|
\ge 2 |f_2|.$ Then $R_{ \sigma } (f_1 + f_2) \le R_{ \sigma }
(f_1) + \frac{2\pi}3. $ Indeed, the value $R_{ \sigma }(f_1 + \varepsilon
f_2)$ cannot change more than by $\frac{2\pi}3$,
as $\varepsilon $ ranges over the segment $[0,1].$

To complete the proof of Proposition \ref{p6.3}, we need to consider
three cases. Let $ \nu = \text {ord}_0 Y, \ \mu = \text {ord}_0 J,
\ f(\ph) = Y ( \psi e^{2\pi i \ph}), \ g (\ph) = \left( J \frac
{\log }{2\pi i}\right) ( \psi e^{2\pi i \ph}).$ Note that $ \mu \le
n^4.$

Case 1: $ \nu < \mu.$ Then, for $\psi$ small, $2 |g| \le |f|.$ By
the previous remark, applied to $f_1 = f, \ f_2 = g, $ we get
$$\left| R_{\Gamma_3} (I)\right| \le \pi (4 \nu + 1)\le \pi (4 n^4 +1).$$

Case 2:  $ \nu = \mu.$ Then, for $\psi$ small, $2 |f| \le |g|,$
because of the logarithmic factor in $g.$ In the same way as
before,  we get
$$
\left| R_{\Gamma_3} (I)\right| \le \pi (4 \mu + 1)\le \pi (4 n^4
+1).
$$

Case 3: $ \nu > \mu.$ In the same way, as in Case 2, we get (\ref{6.6}).
\end{proof}



\subsection{Proof of Theorem A2 in case 1 (endpoints of the
interval considered are finite)}

\begin{proof} It is sufficient to prove Lemma \ref{l6.2}.
We prove a stronger statement
\begin{equation}
N(I,\Pi'_\psi ):= \# \{t \in \Pi'_\psi \mid I(t) =0\} < \frac
{1}{2}(1 -\log c')A^{4700}   \label{6.8}
\end{equation}
By the argument principle
\begin{equation}
2\pi N(I,\Pi'_\psi ) \le V(\Gamma_1) + \mid R_{\Gamma'_2}(I)\mid +
\mid R_{\Gamma_3}(I)\mid  + \mid R_{\Gamma'_4}(I)\mid   \label{6.9}
\end{equation}
The first term in the r.h.s is estimated in (\ref{6.3}). The
second and the fourth terms are estimated from above in
(\ref{6.5a}). The third term is estimated in (\ref{6.6}).
Altogether this proves (\ref{6.8}), hence, Lemma 4.2 and implies a
stronger version of (\ref{6.1}):
$$
\# \{t \in \Pi'_\psi \mid I(t) =0\} < \frac {1}{2}(1 -\log
c')A^{4700}.
$$

This proves Theorem A2 in case 1.
\end{proof}

\subsection{Proof of Theorem B2}
\def\G{\Gamma}
\medskip
{\bf Theorem B2.} {\it For any real ultra-Morse polynomial $H$, any
family $\G $ of real ovals of $H$, and any $l,$ let $\Pi (a)$ and
$\Pi (b)$, $D(l,a)$ and $D(l,b)$ be the same domains,
as in (\ref{1.2b}). Let
$I$ be the analytic extension to $W$ of the integral (1.1)
over the ovals of the family $\G :$
$$
\ \ \ \ \int_{\gamma_t} \omega =I(t),\ \gamma_t \in \Gamma .
$$
Then the number of zeros of $I$ in $D(l,a)$ and $D(l,b)$
(denoted by $N(l,H)$), is no greater than}
$$
N(l,H) \le (1 - \log c'(H))e^{4700\frac {n^4}{c''}+\frac
{481l}{c''}}.
$$
\medskip

\begin{proof} We will prove the theorem for the case when $a =
a(t_0)$ is a logarithmic branch point of the integral $I$ at the
left end of the segment $\sigma (t_0).$ The case of the right end
is treated in the same way. The case when $a(t_0)$ is a critical
value of $H$ which is not a singular point of the integral $I,$ is
even more elementary. In this case the integral is univalent in a
small \nbd of $a,$ the number of zeros to be estimated does not
depend on $l \ge 1,$ and the estimate follows from TheoremA2.

Let for simplicity $D(l) =D(l,a).$ For any $\psi \in (0, \nu
)$ consider the set $\Pi'_\psi \subset W,$ see (\ref{6.4a}). Let
$$
\Pi'_{\psi ,l} = \{ re^{i\ph } + a(\G ) \in W \mid re^{\frac {i\ph
}{l}} +a(\G ) \in \Pi'_\psi \} .
$$
Let $\G_{1,l}$, $\G'_{2,l}$, $\G_{3,l}$, $\G'_{4,l}$ be the curves 
defined by the relations:
$$
\partial \Pi'_{\psi ,l} = \G_{1,l}\G'_{2,l}\G_{3,l}\G'_{4,l};
$$
$$
\pi \G_{j,l} = \pi \G_j, \ j = 1; 3;
$$
$$
\pi \G'_{j,l} = \pi \G'_j, \ j = 2; 4.
$$
Let $R_\G (f)$ and $V_\G (f)$ be the same as in 4.1. Then, by the
argument principle
\begin{equation}
2\pi N(l,H) \le V_{\G_{1,l}}(I) +\mid R_{\G'_{2,l}}(I)\mid + \mid
R_{\G_{3,l}}(I)\mid + \mid R_{\G'_{4,l}}\mid.    \label{6.9a}
\end{equation}
The four terms in the right hand side are estimated in a similar
way as the corresponding terms in (\ref{6.9}). The last three terms are
in fact already estimated:
\begin{equation}
\mid R_{\G'_{j,l}}(I)\mid \le \pi (2N + 1), \ j = 2,4,  \label{6.10}
\end{equation}
where $N$ is the same as in (\ref{6.4});
\begin{equation}
\mid R_{\Gamma_{3,l}'}(I)\mid \le \pi (4n^4l + 1).   \label{6.11}
\end{equation}
\end{proof}
\begin{proposition} \label{p6.4} Inequality (\ref{6.10}) holds.
\end{proposition}
\begin{proof} The proposition is proved in the very same way as
Proposition \ref{p6.2} with the only difference: (\ref{6.5b})
should be
replaced by
$$
\text { Im }I \mid_{\G_{j,l}} =\pm ll_0J\circ \pi \mid_{\G_0}.
$$
The factor $l$ in the r.h.s. does not change the number of zeros.
\end{proof}

Inequality (\ref{6.11}) is proved in the same way as (\ref{6.6})
with the only
difference that the increment of the argument of $t$ along
$\G_{3,l}$ is now $4\pi l.$

\begin{proposition} \label{p6.5} Let $A = e^{\frac {n^4}{c''}}.$ Then
\begin{equation}
V_{\G_{1,l}} (I) \le (1 - \log c'(H))A^{4700}e^{\frac {481l}{c''}}
\label{6.12}
\end{equation}
\end{proposition}
\begin{remark} Inequalities (\ref{6.10}), (\ref{6.11}),
(\ref{6.12}) together prove Theorem B2.
\end{remark}
\begin{proof} {\bf of proposition \ref{p6.5}.}
The proof follows the same lines
as that of Lemma \ref{l6.3}. We will estimate
the variation of argument
under consideration making use of Theorem \ref{t6.2}.
For this we need
first to choose the curve $\G $ and domains $U'', U', U.$ Let
$$
\G = \G_{1,l} = \{ a + \nu e^{i\ph }\mid \ph \in [-2\pi l, 2\pi
l]\} .
$$
Take the same $\e $ as in (\ref{6.3'}).
For any set $A \subset W$ take
$A^\e $ to be the $\e $-\nbd of $A$ in the Euclidean metric of $W,$
and $\overline {A^\e }$ be the minimal simply connected domain
that contains $A^\e .$ Let $K$ be the same as in the Main Lemma.
Take
$$
U'' = \overline {{(K \cup \G )}^\e }, \ U' = \overline {{(K \cup
\G )}^{2\e }}, \ U = \overline {{(K \cup \G )}^{3\e }}.
$$
Note that for any point $p \in K \cup \G ,$ the $6\e $-\nbd of $p$ in
$W$ is bijectively projected to a $6\e $-disk in $\Bbb C.$ Hence,
the gap condition (\ref{4}) holds for $\G , U'', U', U$ so chosen.

Note that $K \cap \G = \G_1 \ne \emptyset .$ Hence, the set $K
\cup \G ,$ as well as $U'', U', U$ is path connected. Then we
have:
$$
\text { diam } K \le 36n^2
$$
by (\ref{1.6}),
$$
\text { diam }_{int} (K \cup \G ) \le 36n^2 + 4\pi l\nu := D_1,
$$
$$
\text { diam }_{int} U'' \le D_1 + 2\e ,
$$
$$
\text { diam }_{int} U' \le D_1 + 4\e .
$$
Hence, diameter condition (\ref{5}) holds with
$$
D_2 = 36n^2 + 16l\nu = 36n^2 + \frac {4lc''}{n^2},
$$
or with $D = 36n^2 +  \frac {4l}{n^2} \geq D_2$ because $c'' \le 1.$

Let us now estimate from above the Bernstein index $B_1 =
B_{U'',U}(I).$ Let $U_0$ be the domain denoted by $U$ in the Main
Lemma.

Then $K\subset U'', \ U_0 \subset U.$ Let $B_0 =B_{K,U_0}(I)$ be
the Bernstein index estimated in the Main Lemma. By (\ref{1.8}),
$$
B_0 < (1 - \log c')A^2, \ A = e^{\frac {n^4}{c''}}.
$$
\end{proof}
\begin{proposition} \label{p6.6}
$$
B_1 \le B_0 + \log (4l+ 1).
$$
\end{proposition}
\begin{proof} By definition,
$$
B_1 =\log \frac {M_1}{m_1},\ B_0=\log\frac{M_0}{m},
$$
where $M_1 = \max_{\baru} |I|, \ m_1 = \max_{\overline{U''}}|I|,
\ M_0=\max_{\overline{U_0}}|I|, m=\max_K|I|.$
Note that $K \subset U''$, hence, $m \le m_1. $

On the other hand, let
$$
M_J =\max|J|\ \text{on the closure of}\ D_{\nu+3\e}(a).
$$
By definition, $\Gamma^{3\e}\subset U_0$,
$\pi\Gamma^{3\e}\subset D_{\nu+3\e}(a)$, $U$ is the minimal
simply connected domain containing $U_0\cup\Gamma^{3\e}$.
By the Picard-Lefschetz theorem and Lemma \ref{l2.4}
$$
M_1 \le M_0 +|l_0|lM_J, \ |l_0| \le 2.
$$
Let us estimate the integral $J$ from above.
Over each point of $\partial D_{\nu+3\e}(a)$ there are two
points of $\partial(\Gamma_1^{3\e})\subset\overline U_0$, where
$\Gamma_1=\Gamma_a$ is the same,
as in (\ref{6.1}).
The difference of the values of $I$ at the two latter
points is equal to
$\pm l_0J$, $0<|l_0|\leq2$. Therefore,
$$
M_J \le 2M_0. \ \text{Hence,}
$$
$$
M_1 \le M_0(4l +1),
$$
$$
B_1 = \log \frac {M_1}{m_1} \leq \log \frac {M_0(4l+1)}{m} = B_0 + \log
(4l + 1).
$$
\end{proof}

Let us now estimate from above other geometric characteristics
used in Theorem \ref{t6.2}, namely, the length and
total curvature of $\G.$ We have:
$$
|\G | = 4\pi l\nu<4 l n^{-2}; \ |\kappa (\G)| \le 4\pi l.
$$
We can now apply Theorem \ref{t6.2}:
$$
V_{\Gamma_{1,l}}(I) \le c_{n,l}e^{\frac {5}{c''}(36n^2+4 l n^{-2})n^2\cdot 24}
$$
where
$$
c_{n,l} = (B_0 + \log (1 + 4l))(\frac{4 l n^{-2}}{\e} +4\pi l + 1).
$$
By the Main Lemma,
$$
B_0 < (1 - \log c'(H))A^2,
$$
where $A = e^{\frac {n^4}{c''}}.$ Elementary estimates imply:
$$
c_{n,l} \le ((1 - \log c'(H))A^3e^l.
$$
This implies (\ref{6.12}).

Together, inequalities (\ref{6.10}) - (\ref{6.12}) imply Theorem B2.

\subsection{Proof of Theorem A2 in Case 2 (near an infinite
endpoint)}

Here we prove Theorem A2 for a segment with one endpoint
(say, $b$) infinity
(statement 2 mentioned at the beginning of the section).

\begin{proposition} \label{p6.7} 
The integral $I$ has
an algebraic branching point at infinity of order $n+1.$
\end{proposition}

\begin{proof} {\bf of Proposition \ref{p6.7}.} 
Let $S_R$ be the circle $|t|=R,$ $R\geq3$, $\Gamma_R$ be the
$(n+1)$ sheet cover of $S_R$ with the base point $-R.$ Consider
the real ovals $\gamma_t$ extended for $t \in W$. For any arc
$\Gamma'\subset\Gamma_R$ going from $-R$ to
$t_{\varphi}=-Re^{i\varphi}$ let $[\Delta_{\Gamma'}]$ be the class
of all the covering homotopy maps $\{ H=-R\}\to\{
H=t_{\varphi}\}$. Let $h$ be the highest homogeneous part of $H$.
If $H=h$, then for any $R$ the class $[\Delta_{\Gamma'}]$ contains
the simple rotation:
$$
R_0: (x,y) \mapsto (e^{\frac
{i\varphi}{n+1}}x,e^{\frac{i\varphi}{n+1}}y)
$$
In the general case, for $R$ large enough the class
$[\Delta_{\Gamma'}]$ contains a map $\Delta_{\Gamma'}$ close to
the rotation. Let us prove this statement. To do this, consider
the extension of the foliation $H=const$ by complex level curves
of $H$ to the projective plane $\Bbb P^2$ obtained by pasting the
infinity line to the coordinate plane $\Bbb C^2$. The foliations
$H=const$ and $h=const$ are topologically equivalent near
infinity. More precisely, for any $r>0$ large enough there exists
a homeomorphism $\Phi$ of the complement $\Bbb P^2\setminus D_r$
($D_r$ is the ball of radius $r$ centered at 0) onto a domain in
$\Bbb P^2$ that  preserves the infinity line such that
$h\circ\Phi=H$. This follows from the statements that the
singularities of these foliations at infinity are the same and of
the same topological type (nodes), and the holonomy mappings
corresponding to circuits around these singularities in the
infinity line are rotations $t\mapsto e^{\frac{2\pi i}{n+1}}t$ in
the transversal coordinate $t=H^{\frac1{n+1}}$. The last statement follows from
the fact that for a generic $C\in \Bbb C$
$$
H(x,y)|_{x=Cy}={(\widetilde Cx)}^{n+1}(1+o(1)),\
\text{as}\ x \to \infty, \ \widetilde C\ne 0.
$$
The homeomorphism $\Phi$ is close to identity near infinity. For
any $r>0$ there exists a $T(r)>0$ such that for any $t$,
$|t|>T(r)$, \ $S_t\cap D_r=\emptyset$. The map
$\Delta_{\Gamma^{\prime}}$ we are looking for is obtained from the
map $R_0$ corresponding to $h$ by conjugation by the homeomorphism
$\Phi$. By construction, its $n+1$- iterate is identity.
\end{proof}

\begin{proof} {\bf of Theorem A2 near infinity.}
Let $V$ be the Riemann
surface of the integral $I.$ Let $\G \subset V$ be the degree $n +
1$ cover of the circle $|t| = 3$ with the base point $t_1 = +3.$
This is a closed curve on $V.$ This curve is a boundary of a
domain on $V$ that covers a \nbd of infinity. Let us denote this
domain by $V_\infty .$ We will estimate from above
$$
N_\infty = \{ t \in V_\infty \mid I(t) = 0 \} .
$$
This will give an upper estimate to the number of zeros of $I$ on
$\sigma^+= (3,+\infty )$ because $\sigma^+ \subset V_\infty .$ We
will use the argument principle in the form
$$
N_\infty \le \frac {1}{2\pi }V_\G (I).
$$

The variation in the right hand side will be estimated by
Theorem \ref{t6.2}. To apply this theorem we need
to define all the entries like in the previous subsection.

Let $\G = \partial V_\infty .$ Without loss of generality we
consider that $I|_{\Gamma}\neq0$ (one can achieve this
by slight contraction of the circle $|t|=3$).
Let $K$ be the same as in the Main
Lemma. Denote by $U_0$ the set $U$ from that lemma: both $K$ and
$U_0$ are taken projected to the Riemann surface
of the integral $I$. Let $\e $
be the same as in (\ref{6.3'}). One has
$K\supset\Sigma$, see 1.7, hence, $K\supset\Gamma$.

Let
$$
U'' = \overline{K^\e}, \ U' = \overline {K^{2\e }}, \ U =
\overline {K^{3\e }}.
$$
By (\ref{1.6}), the diameter condition (\ref{5}) holds with
$$
D = 36n^2+1.
$$

The gap condition (\ref{4}) for $\G , U'', U', U$ holds as well. The
Bernstein index $B = B_{U'',U}(I)$ may be easily estimated with
the use of the same results that were used in the estimate of $B_0
= B_{K,U_0}(I).$ Indeed,
$$
B = \log \frac {M'}{m'}, \ M' = \max_{\baru}|I|, \ m' =
\max_{\overline{U''}}|I|;
$$
$$
B_0 = \log \frac {M_0}{m}, \ M_0 = \max_{\overline{U_0}}|I|,
\ m = \max_K|I|.
$$
But $K \subset U'';$ hence, $m' \ge m.$ On the other hand,
$$
U = U_0 \cup \G^{3\e }.
$$
Each points of $U$ may be connected to $t_0$ by a path that
satisfies the assumptions of Theorem \ref{t2.1}.
(Let $M_0$ be the constant from the same theorem.)
Hence, by Corollary \ref{cor2},
$$
M' = \max_U\mid I(t)\mid < n^4M_0=M_1'
$$
In the proof of the Main Lemma we used the following inequalities:
$$
M' \le M_1', \ \log \frac {M_1'}{m} < (1 - \log c')A^2.
$$
Hence,
$$
B = \log \frac {M'}{m'} < \log \frac {M_1'}{m} < (1 - \log c')A^2.
$$
This inequality will be substituted in (\ref{3}). Another quantities
from (\ref{3}):
$$
e^{\frac {5D}{\e }} \le A^{4700},
$$
$$
\mid \G \mid \le 6\pi (n + 1),
$$
$$
\mid \kappa (\G )\mid \le 2\pi (n + 1).
$$
Altogether, by Theorem \ref{t6.2}, this implies Theorem A2, Case 2.
\end{proof}

\section{Acknowledgements}
The authors are grateful to L.Gavrilov,
P.Haissinsky and S.Yu.Yakovenko for helpful discussions.

\end{document}